\documentclass[12pt, reqno]{amsart}
\usepackage{amscd,times}
\usepackage{amssymb}
\usepackage{amsmath}
\usepackage[all]{xypic}
\usepackage{graphicx}

\newtheorem{thm}{Theorem}[section]
\newtheorem{lem}[thm]{Lemma}

\newtheorem{cor}[thm]{Corollary}
\newtheorem{pro}[thm]{Proposition}
\newtheorem{Exa}[thm]{Example}
\newcommand{\bbar}{\begin{array}}
\newcommand{\eear}{\end{array}}
\newcommand{\bb}{\begin{equation}}
\newcommand{\eqbb}{\begin{equation}}
\def\ee{\end{equation}}
\def\eqee{\end{equation}}
\def\eea{\end{eqnarray}}
\def\bba{\begin{eqnarray}}

\def\ch{\mbox{\rm ch}}

\def\secteqn{
\let\sectio\section%
\renewcommand{\section}{\sectioneqn\sectio }%
\newcommand{\sectioneqn}{\setcounter{equation}{0}%
 \renewcommand{\theequation}{\arabic{section}.\arabic{equation}}}}\oddsidemargin0.5cm
\begin{document}

\title[Jacobi-Trudi type formula]{Jacobi-Trudi type formula for a class of irreducible representations of  $\frak{gl}(m|n)$ }
  
\author{Nguy\^en Luong Th\'ai B\`inh}
 \address[NLT Binh]{Sai Gon University, Ho Chi Minh City, Vietnam and Institute of Mathematics, Vietnam Academy of Science and Technology, Hanoi, Vietnam}
 \email{nltbinh@sgu.edu.vn}
%

%
%

 \maketitle
 \bibliographystyle{plain}
\begin{abstract}
We prove a determinantal type formula to compute the characters for a class of irreducible representations of the general Lie superalgebra $\mathfrak{gl}(m|n)$ in terms of the characters of the symmetric powers of the fundamental representation and their duals. This formula was conjectured by J. van der Jeugt and E. Moens and was generalized the well-known Jacobi-Trudi formula.
\end{abstract}

 \section{Introduction}

The classical Jacobi-Trudi formula computes Schur symmetric functions in terms of the elementary (resp. complete) symmetric functions. Since these symmetric functions can be realized as irreducible characters of a general linear group, we can interpret the Jacobi-Trudi formula as a formula for computing irreducible character of a general linear group in terms of the characters of symmetric (resp. anti-symmetric) tensor representations. This formula complements the Weyl determinantal formula which computes irreducible characters in terms of the root system.
Although the Jacobi-Trudi formula is well-defined only for partitions, that is, for integral dominant weight with non-negative components, it is well-known that an integral dominant weight can be led to a partition by adding some multiple of the partition $(1,1,\ldots,1)$, which corresponds to the  determinantal representation.

 The aim of this work is to extend this famous formula to the case of the general linear Lie super algebras. According to V.~Kac, irreducible representations of the general linear Lie super-algebra $\mathfrak{gl}(m|n)$ are determined by means of dominant weights. We shall restrict ourselves to those representations with integral dominant weights, see Eq. \eqref{eq:lamda1}. In his foundational papers \cite{Kac1,Kac2,Kac3,Kac4} on Lie superalgebras, Kac raised the problem of determining the formal characters of finite dimentional irreducible representations of Lie superalgebras, and he established in the 70s an analog of Weyl formula to compute irreducible characters corresponding to typical weights. There have been many attempts to establish formulas to computing atypical irreducible characters, see e.g. \cite{Brundan,bar1,bar2,CK,dj, M1,M2,M3,M4,M5,M6,Zhang2,Jeugt1,Jeugt2}. It took twenty years until V. Serganova provides a method to compute atypical irreducible characters, which was subsequently simplified by Brundan \cite{Brundan} and Su-Zhang \cite{Zhang2}.

For some classes of integral dominant weights, Jacobi-Trudi formula has been established, for instance, when the weights correspond to partitions, i.e. the corresponding representation is constructed from the fundamental representation using only tensor products and decomposition into direct sums, see \cite{br1, dj}. However, due to the more complicated nature of the representation category of general linear Lie superalgebra, to extend Jacobi-Trudi formula to characters of the so-called mixed representations one will need to incorporate characters of both symmetric tensor powers and their duals. A conjectural determinantal formula was explained in detail in \cite{M6}. In fact, there was an unsuccessful attempt to prove it in \cite{M4}. This formula can be considered as an analog of the Jacobi-Trudi formula for irreducible characters of the general linear Lie superalgebras.

In this work we prove the above mentioned determinantal formula for the case of  irreducible representation correspond to integral dominant weight $\Lambda$ which has the form \\
 $$\Lambda = (\alpha_1, \alpha_2, \ldots, \alpha_m; -k, -k, \ldots,  -k)$$
 such that $0\leq k\leq m $ and $\alpha_{m-k} \geq 0 \geq \alpha_{m-k+1}$ (Theorem \ref{Mt}). A class of  these weights corresponds to a class of $m$-standard composite partitions (cf. \eqref{eq:42}). In particular, we have also derived the results presented in \cite{bdh}, that is Jacobi-Trudi type formula for character of  irreducible representations of  $\frak{gl}(m|1)$ (Corollary \ref{C:cor2}).

The structure of the paper is as follows. In section 2 we present some background materials on $\mathfrak{gl}(m|n)$. In section 3 we state the main theorem (see Theorem \ref{Mt}) and some corollaries (Corollary \ref{C:cor1} and Corollary \ref{C:cor2}). In section 4 we introduce the notion special weights (cf. \eqref{eq:sp}) and establish the correspondence between these weights and the $(m|n)$-standard composite partitions (Proposition \ref{Pro:2}). In section 5 we give a reduction formula to represent irreducible characters of $\frak{gl}(m|n)$ in terms of irreducible characters of its subalgebras (Theorem \ref{T:ch}). In section 6, we recall the notion of super-symmetric $S$-functions corresponding to a composite partitions and properties their (see Theorem \ref{T:S-function}). The last section gives a proof of the main theorem.

\section{Preliminaries}
This section presents some results on linear Lie super-algebras for the later use. We shall work over the field of complex numbers $\mathbb{C}$. \\

A super vector space is a $\mathbb{Z}_2$-graded vector space $V = V_{\bar{0}}\oplus V_{\bar{1}}$. The vector spaces $ V_{\bar{0}}, V_{\bar{1}}$ are called the even and odd homogeneous components of $V$, their elements are also called homogeneous. A homogeneous element $x \in V_{\bar{0}}$ has degree 0, denoted $deg(x) = \bar{0}$, while $x \in V_{\bar{1}}$ has degree 1, denoted $deg(x) = \bar{1}.$\\

Let $End(V)$ be the space of linear endomorphisms of $V$. Then $End(V) = End_{\bar{0}}(V) \oplus End_{\bar{1}}(V)$, where
\begin{equation}\label{eq:}
End_{\bar{0}}(V) = End(V_{\bar{0}}) \oplus End(V_{\bar{1}})\;\; \mbox{and}\;\; End_{\bar{1}}(V) = Hom(V_{\bar{0}}, V_{\bar{1}}) \oplus Hom(V_{\bar{1}}, V_{\bar{0}}).
\end{equation}
We can equip $End(V)$ with the structure of a Lie superalgebra by defining the Lie bracket  $[-,-]$ by setting
\begin{equation}\label{eq:}
[x,y] = xy - (-1)^{deg(x)deg(y)}yx,
\end{equation}
for all homogeneous elements $x, y$, then extending it linearly to the whole space $End(V)$ with $V = V_{\bar{0}}\oplus V_{\bar{1}}$ and dim$V_{\bar{0}} = m$, dim$ V_{\bar{1}} = n$.
 We use the notation $\frak{gl}(m|n)$ for the $End(V)$ with above Lie superalgebra structure.

\subsection{The Lie superalgebra $\frak{gl}(m|n)$}

In this paper, $\frak{g}$ will denote the Lie superalgebra $\frak{gl}(m|n)$ for fixed positive integers $m, n$. We can realize $\frak{g}$ as the set of $ (m + n) \times (m + n)$ matrices.
Hence 
\begin{equation}\label{eq:}
\frak{g}_{\bar{0}} = \left\{ \left (\begin{array}{cc}A&0\\
0&D \end{array}\right) | A \in M_{m,m}, D \in M_{n,n} \right\} \end{equation}
and
 \begin{equation} \frak{g}_{\bar{1}} = \left\{ \left (\begin{array}{cc}0&B\\
C&0 \end{array}\right) | B \in M_{m,n}, C \in M_{n,m} \right\},
\end{equation}
where $M_{r,t}$ denotes the set of $r\times t$ matrices.\\

The standard basis for $\frak{g}$ consists of the matrices $E_{i,j}, i,j = 1, 2, \ldots, m + n$, with 1 on the entry $(i, j)$ and $0$ elsewhere. The subalgebra $\frak h \subset  \frak g$ spanned by the elements $ E_{j,j}: j = 1, 2, \ldots, m + n$, is called Cartan subalgebra of $\frak{g}$. The dual vector space  $\frak{h}^*$ is called the weight space, it  is spanned by the weights $\{\epsilon_i, \delta_j| i = 1, 2, \ldots, m; j = 1, 2,\ldots, n  \}$, where $\epsilon_i(E_{j,j}) = \delta_{ij}$ and $\delta_j(E_{i,i}) = - \delta_{(m + j)i}$, with $\delta$ is symbol Kronecker .

A weight $\Lambda$ will be denoted as follows:
\begin{equation}\label{eq:lamda1}
\Lambda 
= \sum_{i=1}^m \lambda_i\epsilon_i +  \sum_{j=1}^n \mu_i\delta_i
=:  (\lambda_1, \cdots, \lambda_m; \mu_1, \mu_2, \ldots, \mu_n).
\end{equation}
$\Lambda$ is called integral if its components (the $\lambda_i's, \mu_j's $) are integers. $\Lambda$  is called dominant if  $\lambda_1 \geq \lambda_2 \geq \cdots \geq \lambda_m$ and $\mu_1\geq \mu_2 \geq \ldots \geq \mu_n$.\\

We fix simple root system 
$$\Pi = \{\epsilon_1-\epsilon_2, \cdots, \epsilon_{m-1}-\epsilon_m, \epsilon_m-\delta_1, \delta_1 - \delta_2, \ldots, \delta_{n-1} - \delta_n \}.$$
The set of positive even roots is denoted by
$$\Delta_{0}^+= \{\epsilon_i - \epsilon_j| 1 \leq i <  j \leq m\} \cup \{\delta_i - \delta_j| 1 \leq i <  j \leq n\} ,$$
 and the set of positive odd roots is denoted by
$$\Delta_{1}^+= \{\epsilon_i - \delta_j | 1 \leq i \leq m, 1 \leq j \leq n \}.$$

As usual, we put
$$\rho_0 =   \frac{1}{2} \sum_{\alpha \in \Delta _0^+} \alpha = \frac{1} {2}(m-1, m-3, \ldots, 1-m; n - 1, n - 3, \ldots, 1 - n),$$
$$  \rho_1 = \frac{1}{2}\sum_{\alpha \in \Delta_{1}^+}\alpha=\frac{1} {2}(n, n, \ldots, n ; -m, -m, \ldots, -m ),$$ 
 $$
 \quad \rho = (m, \ldots, 2, 1; -1, -2, \ldots, -n).
$$

There is a symmetric bilinear form $(\;,\;)$ on $\frak h^*$ is defined by
$$(\epsilon_i, \epsilon_j) = \delta_{ij}, (\epsilon_i, \delta_j) = 0, (\delta_i, \delta_j) = -\delta_{ij}.$$
The Weyl group of $\frak g $ is the Weyl group $ W$ of $\frak g_0$, hence it identified with the product of the symmetric groups $ S_m\times S_n $. For $ w \in W$, we denote by $\epsilon(w)$ its signature.\\
\subsection{Typical and atypical weights} 
Let $\Lambda = (\lambda_1, \cdots, \lambda_m; \mu_1, \mu_2, \ldots, \mu_n)$ be an integral dominant weight.
A positive odd root $\epsilon_i - \delta_j $, with $i = 1, 2, \ldots, m$ and $j = 1, 2, \ldots, n$, is called an atypical root of $\Lambda$ if\\
\begin{align}
(\Lambda + \rho, \epsilon_i - \delta_j) = 0.
\end{align}
Explicitly, this condition reads: $  \lambda_i +m + 1 - i =  - \mu_j + j$.
Denote by $\Gamma_\Lambda$ the set of atypical roots of $\Lambda$:\\
$$\Gamma_\Lambda = \{\epsilon_i - \delta_j| (\Lambda + \rho, \epsilon_i - \delta_j) = 0 \}.$$
 A weight $\Lambda$ is called typical if $\#\Gamma_\Lambda =0$ and atypical if  $\#\Gamma_\Lambda = r \geq 1$ (in this case $\Lambda$ is also called $r$-fold atypical weight).

\subsection{Kac modules}
For every integral dominant weight $\Lambda$, we denote by $V^0(\Lambda)$ the finite dimensional irreducible $ \frak{g}_{\bar{0}}$-module with highest weight $\Lambda$. $V^0(\Lambda)$ is a $(\frak{g}_{\bar{0}} \oplus \frak{g}_{+1})$- module with $\frak{g}_{+1}$ acting by $0$, where $\frak{g}_{+1}$ is the set of matries of the form $\left (\begin{array}{cc}0&B\\
0&0 \end{array}\right)$. Set
$$ 
\bar{V}(\Lambda) := \text{\rm Ind}_{\frak{g}_{\bar{0}} \oplus \frak{g}_{+1}}^{\frak{g}}V^0(\Lambda).
 $$
Then $\bar{V}(\Lambda)$ contains unique maximal submodule $M(\Lambda)$. So that, $\bar{V}(\Lambda)/M(\Lambda)$ is an irreducible module. Put
$$ V(\Lambda) := \bar{V}(\Lambda)/M(\Lambda).  $$
Then $V(\Lambda)$ is an irreducible module with highest weight $\Lambda$. It is called Verma module or Kac module \cite{Kac1}. 

\subsection{Characters of $\frak{gl}(m|n)$}
Let $V(\Lambda)$ be an irreducible representation with highest weight $\Lambda$ of $\frak g$. Such a representation is $\frak h$-diagonalizable with weight decomposition $ V(\Lambda) = \bigoplus_{\mu}V_\mu$, where $V_\mu = \{ v \in V | h v = \mu(h) v \;\; \mbox{for all}\;\; h \in \frak{h}\}$. The character of  $V(\Lambda)$ is defined to be the formal sum
$$\ch \; V = \sum_{\mu}(\dim V_\mu) e^{\mu},$$
 where $e^{\mu}$ ($\mu \in \frak h^*$) are the formal exponential functions.

 \section{The Main theorem}
In this section we will state the main theorem, the Jacobi-Trudi type formular to compute characters of a class of irreducible representations of a general linear Lie superalgebra $\mathfrak{gl}(m|n)$. This class consists of irreducible representations correspond to integral dominant weights $\Lambda$ of the form 
$$ \Lambda = (\alpha_1, \alpha_2, \ldots, \alpha_m; -k, -k, \ldots,  -k),$$
  with $0\leq k\leq m $ and $\alpha_{m-k} \geq 0 \geq \alpha_{m-k+1}$.

Set $x_i= e^{\epsilon_i}$ for $i =1,2,\ldots,m$ and $y_j= e^{\delta_j}$ for $j = 1,2,\ldots, n$. It is well known that the character of the $r$-th the symmetric power of the fundamental representation of  $\mathfrak{gl}(m|n)$ is equal to 
$$ h_r(x/y) := \sum_{k = 0}^r h_k(x)e_{r-k}(y),$$
where $e_k(y)$ (resp. $h_k(x)$) is the elementary (resp. complete) symmetric function on the variables $ y = (y_1, y_2, \ldots, y_n)$ (reps. $ x= (x_1, x_2, \ldots, x_m) $). 
It follows that the character of the dual of the $r$-th symmetric power of the fundamental representation is $ \dot{h_r}(x/y) $ with $\dot{h_r}(x/y): = h_r(\bar{x}/\bar{y})$ where $\bar{x}= (x_1^{-1},x_2^{-1},\ldots,x_m^{-1}), \bar{y} = (y_1^{-1},y_2^{-1},\ldots,y_m^{-1})$.

\begin{thm}\label{Mt}
Let  $ \Lambda = (\alpha_1, \alpha_2, \ldots, \alpha_m; -k, -k, \ldots,  -k)$ be an integral dominant weight of a general linear Lie superalgebra $\mathfrak{gl}(m|n)$ with $0\leq k\leq m $ and $\alpha_{m-k} \geq 0 \geq \alpha_{m-k+1}$. Then
$$  {\rm ch} V(\Lambda) = det\left(
\begin{tabular}{c|c}
$\dot{h}_{n-\alpha_{m-t+1}+s-t}(x/y)$&$h_{\alpha_j-s-j+1}(x/y)$\\
\hline
$\dot{h}_{n-\alpha_{m-t+1}-i-t+1}(x/y)$&$h_{\alpha_j +i -j}(x/y)$
\end{tabular}\right),$$
where the indices $i, j, s$ resp. $t$ run from top to bottom, from left to right, from bottom to top resp. from right to left with $i,j = 1,2,\ldots, m-k$ and $s,t=1,2,\ldots, k$.
\end{thm}
Example: For $\Lambda = (3,2,-1;-1-1)$, an integral dominant weight of the linear Lie superalgebra $\frak{gl}(3|2)$, we have
$$  {\rm ch} V(\Lambda) = \left|\begin{array}{ccc}
\dot{h_3}(x/y)&h_2(x/y)&h_0(x/y)\\
\dot{h_2}(x/y)&h_3(x/y)&h_1(x/y)\\
\dot{h_1}(x/y)&h_4(x/y)&h_2(x/y)
\end{array}
\right|.$$

Consider the weight $\sigma = (1,1,\ldots, 1; -1, -1, \ldots, -1)$ which corresponds to the super-determinantal representation. Let $ \Lambda = (\lambda_1, \lambda_2, \ldots, \lambda_m; \beta,\beta,\ldots , \beta )$ be an integral dominant weight. Then there is a unique integer $j$ such that 
$$ \Lambda + j\sigma = \Lambda_0,  $$
where $\Lambda_0 = (\alpha_1, \alpha_2, \ldots, \alpha_m; -k, -k, \ldots,  -k)$ with $0\leq k\leq m $ and such that $\alpha_{m-k} \geq 0 \geq \alpha_{m-k+1}$ i.e, $\Lambda_0 $ satisfies the condition of Theorem \ref{Mt}(see Proposition \ref{Pro:1}). We present an algorithm to find $\Lambda_0$, i.e, to find $j$ and $k$. \\
First, we see that, if $ \Lambda = (\lambda_1, \lambda_2, \ldots, \lambda_m; \beta,\beta,\ldots , \beta )$ then
$$ \Lambda+ \beta \sigma = (\lambda_1 + \beta, \lambda_2+ \beta, \ldots, \lambda_m+\beta; 0,0, \ldots , 0 ). $$ 
Without  lost of generality, we can consider $ \Lambda = (\lambda_1, \lambda_2, \ldots, \lambda_m; 0,0,\ldots , 0 )$ instead of $ \Lambda = (\lambda_1, \lambda_2, \ldots, \lambda_m; \beta,\beta,\ldots , \beta )$. Now, we will show an algorithm for finding $j$ and $k$ :
$$ \Lambda + j\sigma =  (\alpha_1, \alpha_2, \ldots, \alpha_m; -k, -k, \ldots,  -k),  $$
 such that $\alpha_{m-k} \geq 0 \geq \alpha_{m-k+1}$. \\

Step 1: If $\lambda_m \geq 0 $ then $j=0$, $k= 0$. If $\lambda_m < 0 $, we move on to Step 2.\\
Step 2: If $\lambda_{m-1}+1 \geq 0 $ then $j=1$, $k= 1$. If $\lambda_{m-1}+1 < 0 $, we move on to Step 3.\\
Step 3: If $\lambda_{m-2}+2 \geq 0 $ then $j=2$, $k= 2$.  If $\lambda_{m-2}+2 < 0 $, we move on to Step 4.\\
$\ldots$\\
Step $m$:  If $\lambda_{1}+ m-1 \geq 0 $ then $j=m-1$, $k= m-1$.  If $\lambda_{1}+ m-1 < 0$ then $j=m$, $k= m$.\\
Thus, after no more than $m$ steps we find $j$ and $k$  that satisfy the requirement.

\begin{cor}\label{C:cor1}
Let $\Lambda = (\lambda_1, \lambda_2, \ldots, \lambda_m; \beta,\beta,\ldots , \beta )$ be an integral dominant weight of $\mathfrak{gl}(m|n)$. Then there is a unique integer $j$ such that
$$ \Lambda + j\sigma = (\alpha_1, \alpha_2, \ldots, \alpha_m; -k, -k, \ldots,  -k),  $$
with $0\leq k\leq m $ and $\alpha_{m-k} \geq 0 \geq \alpha_{m-k+1}$. And

\begin{equation}\label{eq:Jt}
{\rm ch} V(\Lambda) =  \left(\frac{\prod_{i=1}^mx_i}{\prod_{j=1}^ny_j}\right)^j det\left(
\begin{tabular}{c|c}
$\dot{h}_{n-\alpha_{m-t+1}+s-t}(x/y)$&$h_{\alpha_j-s-j+1}(x/y)$\\
\hline
$\dot{h}_{n-\alpha_{m-t+1}-i-t+1}(x/y)$&$h_{\alpha_j +i -j}(x/y)$
\end{tabular}\right),
\end{equation}
where the indices $i, j, s$ resp. $t$ run from top to bottom, from left to right, from bottom to top resp. from right to left with $i,j = 1,2,\ldots, m-k$ and $s,t=1,2,\ldots, k$.
\end{cor}

\begin{proof}
This follows from the Theorem \ref{Mt} and from the following formula, which is well-known:
$${\rm ch} V(\Lambda + j\sigma) = (e^{\sigma})^j{\rm ch} V(\Lambda).$$
\end{proof}

 We now apply the above result to the case $\mathfrak{gl}(m|1)$. The corollary below is proven by direct computation in \cite{bdh}.
\begin{cor}\label{C:cor2}
Let $V$ be an arbitrary irreducible representation of a linear Lie superalgebra $\mathfrak{gl}(m|1)$. Then  ${\rm ch} V$ is the product of the power of $y^{-1}\prod_{i=1}^mx_i$ and the Jacobi-Trudi type formular.
\end{cor}

\begin{proof}
Any integral dominant weight of  $\mathfrak{gl}(m|1)$ satisfies the condition of Corollary \ref{C:cor1}.
\end{proof}

\section{Special weights and composite partitions}
 \subsection{Special weights} 
An integral dominant weight 
\begin{align}\label{eq:sp}
  \Lambda = (\alpha_1, \alpha_2, \ldots, \alpha_m; \beta_{1}:= -k , \beta_2, \ldots,  \beta_ n),
\end{align}
  with $ 0 \leq k \leq m$ and $\alpha_{m-k}\geq 0 \geq \alpha_{m-k+1}$ , is called a \textit{special weight}.
 We denote by $P$ the set of all special weights. And for each integer $k$, $0 \leq k \leq m$, set
\begin{equation}\label{eq:15}
P_k = \{  \Lambda = (\alpha_1, \alpha_2, \ldots, \alpha_m; \beta_{1}:= -k , \beta_2, \ldots,  \beta_ n)  | \alpha_{m-k}\geq 0 \geq \alpha_{m-k+1}     \}.
\end{equation} 
It's easy to see that
\begin{equation}\label{eq:a}
 P = \bigcup_{k=0}^mP_k .
\end{equation} 

  We will see that an arbitrary integral dominant weight of $\frak{gl}(m|n)$ is different from a weight in $P$ by an integer multiple of the weight $\sigma:=(1,\ldots,1;-1, \ldots, -1)$.
\begin{pro}\label{Pro:1}
Let $\lambda $ be an integral dominant weight. Then, there is unique integer $j$ such that $\Lambda :=\lambda + j\sigma $ has the following form:
$ \Lambda = (\alpha_1, \alpha_2, \ldots, \alpha_m; -k, \beta_2, \ldots,  \beta_ n)$ with $ 0 \leq k \leq m$ and  $\alpha_{m-k}\geq 0 \geq \alpha_{m-k+1}$.
\end{pro}
 
\begin{proof}
Existence of $j$. We use induction on $m$. Let's consider the case $ m = 1$. Consider a weight
$(\alpha_1; \beta_1, \beta_2, \ldots,  \beta_ n)$. Then by adding a mutiple of the weight $\sigma$: $j = \beta_1$ if $\alpha_1 +  \beta_1 \geq 0$ and $j = \beta_1 + 1$ if $\alpha_1 +  \beta_1 < 0$, we have
$$\left[\begin{array}{ll} (\alpha_1; \beta_1, \beta_2, \ldots,  \beta_ n) + \beta_1\sigma & \text{ if } \alpha_1 +  \beta_1 \geq 0,\\
 (\alpha_1; \beta_1, \beta_2, \ldots,  \beta_ n) + ( \beta_1 + 1)\sigma & \text{ if }   \alpha_1 +  \beta_1 < 0.\end{array}\right.$$
this is equivalent to
$$
\left[\begin{array}{ll} (\alpha_1 + \beta_1; 0, \beta_2  - \beta_1, \ldots,  \beta_ n - \beta_1)  & \text{ if } \alpha_1 +  \beta_1 \geq 0,\\
 (\alpha_1 +  \beta_1 + 1; -1, \beta_2 - \beta_1 -1 , \ldots,  \beta_ n -  \beta_1 -1)  & \text{ if } \alpha_1 +  \beta_1 < 0.\end{array}\right.$$

 Suppose this holds for $m$. Consider a weight
$$  (\lambda_1,\lambda_2, \ldots, \lambda_{m + 1}; \nu_1, \nu_{2},\ldots,  \nu_n), $$
using the induction hypothesis on the weight
 $$(\lambda_2,\ldots, \lambda_{m + 1}; \nu_1, \nu_{2}, \ldots,  \nu_n),$$
we can bring it to the form
$$  (\alpha_1, \alpha_2, \ldots, \alpha_m, \alpha_{m + 1}; -k, \beta_2, \ldots,  \beta_ n), $$
with $0\leq k\leq m$ (by adding a multiple of $\sigma$), such that
$$\alpha_{(m+1)-k}\geq 0 \geq \alpha_{(m+1)-k+1}.$$
Note the shift of indices and the condition $\alpha_1\geq 0$ is not imposed  (when $k = m$) .\\

Thus, if in this new weight we have $k < m $, then it automatically satisfies the requirement. Similarly, if in this new weight we have $k = m$ and $\alpha_1\geq 0 $ then it also satisfies the requirement. 
It remains the case $k = m$ and $\alpha_1 < 0 $. Then adding $\sigma$ to this weight we get the weight 
$$ (\alpha_1 + 1, \alpha_2 + 1,\ldots, \alpha_{m + 1} + 1; -( m + 1) ,\ldots, \beta_{n}- 1)$$ 
such that $ 0 \geq \alpha_1 + 1$.\\

Finally, we prove the uniqueness assertion. Let $\lambda$ be a integral dominant weight. Assume $j, j'$ are integers such that

$$
 \begin{cases}
   \lambda + j\sigma = (\lambda_1, \lambda_2,\ldots, \lambda_{m };- k, \nu_2,\ldots, \nu_{n} ), &\mbox{with}\;\; \lambda_{m-k}\geq 0 \geq \lambda_{m - k +1}\\ 
     \lambda + j'\sigma = (\alpha_1, \alpha_2, \ldots, \alpha_m; -k', \beta_2, \ldots,  \beta_ n), &\mbox{with}\;\; \alpha_{m-k'}\geq 0 \geq \alpha_{m - k' +1}.
\end{cases} 
 $$
We need to show that  $ j = j' $. Assume the contrary $i \ne j$, then we can say $j > j'$. We have

\begin{align*}
(j - j')\sigma &= \lambda + j\sigma -(\lambda + j'\sigma)\\ 
&=  (\lambda_1- \alpha_1, \lambda_2 - \alpha_2,\ldots, \lambda_{m }- \alpha_m; - (k - k'), \nu_2- \beta_2,\ldots, \nu_{n}-\beta_{n}).
\end{align*} 
Then 
$$ j -  j' = k - k' := t > 0.$$

On the other hand,
\begin{align*}
\lambda + j\sigma &= \lambda + j'\sigma + t\sigma\\ 
& = (\alpha_1 + t, \alpha_2 + t,\ldots,  \alpha_m + t; - ( k'+ t), \beta_2 - t,\ldots, \beta_{n}- t).
\end{align*}
Thus
$$ \alpha_{m - k} + t \geq 0 \geq \alpha_{m-k +1} + t .$$
This implies
$$  0 \geq \alpha_{m-(k' + t) +1} + t \geq \alpha_{m - k' } + t \geq  t  > 0,$$
 which is a contradiction. Thus we conclude that $j = j'$.

\end{proof}

\subsection{Composite partitions}
%
%

Let $\nu, \mu$ be two partitions. We shall refer to $\bar{\nu};\mu$ as a \textit{composite partition}.\\
 A composite partition is called an \textit{$m$-standard composite partition} if  
\begin{equation}\label{eq:42}
l(\mu) + l(\nu) \leq m.
\end{equation}

A composite partition $\bar{\nu};\mu$ is said to be an \textit{$(m|n)$- standard} if  there exist $0 \leq j \leq n$ and $0 \leq l \leq m$ such that
$$
\begin{cases}
 &\mu '_{j+1} + \nu'_{n-j+1} \leq m \\
&\mu_{m-l+1} + \nu_{l+1} \leq n ,
\end{cases} $$
where $\mu'$ (resp. $\nu'$) is conjugate partition of $\mu$ (resp. $\nu$).\\

For each $ 0 \leq k \leq m$
let $Q_k$ be the subset of $(m |n)$-standard composite partitions $\bar{\nu};\mu$, for which $\mu'_1 \leq m - k$ and $\nu'_n =k$:
\begin{equation}\label{eq:}
 Q_k = \{ \bar{\nu};\mu | \mu'_1 \leq m - k , \nu'_n =k\}.
\end{equation} 
Put
\begin{equation}\label{eq:}
 Q = \bigcup_{k=0}^mQ_k .
\end{equation} 
We define a map $\varphi_{k}:P_{k}\to Q_k$ as follows. For a $\Lambda = (\alpha_1, \alpha_2, \ldots, \alpha_m; -k, \beta_2, \ldots,  \beta_ n) \in P_k$ with $ 0 \leq k \leq m$. 
 $\varphi_k(\Lambda)$ is the composite partition $\bar{\nu};\mu$ where the partition $\mu = (\alpha_1, \alpha_2, \ldots, \alpha_{m - k})$ and 
the partition $\nu$ is given by $\nu_1 = n-\alpha_m, \nu_2= n-\alpha_{m-1}\ldots, \nu_k = n- \alpha_{m-k+1}$ while  $\nu_{k+1}, \nu_{k+2}, \ldots$ are uniquely determined by 
$\nu'_n =k, \nu'_{n-1}= -\beta_2, \ldots, \nu'_1 = - \beta_n .$
Notice that $\mu'_1 \leq m - k$ and $\nu_n'=k$. Thus, if $\Lambda\in P_k$ then  
$ \varphi_k(\Lambda) = \bar{\nu};\mu \in Q_k.$  

\begin{pro}\label{Pro:2}
 $\varphi_k$ defines a bijection between $P_k$  and $Q_k$ for each $0\leq k\leq m$. Consequently we have a bijective map $\varphi$ from $P$ to $Q$ such that $\varphi_{|P_k} = \varphi_k$.
\end{pro}
\begin{proof}[Proof]
It is easy to see that $\varphi_k$ is injective. \\
For the surjectively, let $ \bar{\nu};\mu \in Q_k$ be an $(m |n)$-standard composite partition. Then $\nu, \mu$ have the following form
$$
\begin{cases}
    \mu &= (\mu_1, \mu_2, \ldots, \mu_{m - k}),   \\ 
    \nu &= (\nu_1, \nu_2, \ldots) \quad \mbox{and}\; \nu'_n =k.
\end{cases} 
$$
We set
$$ \Lambda = (\mu_1, \mu_2, \ldots, \mu_{m - k}, n-\nu_k, \ldots, n - \nu_2, n-\nu_1; -k=-\nu'_n, \ldots, -\nu'_2, -\nu'_1).
 $$
This weight have $\mu_{m - k}\geq 0 \geq n-\nu_k$ so $\Lambda \in P_k$.
It is obviously the preimage of $\bar{\nu};\mu$ under $\varphi_k$. Thus $\varphi_k$ and $\varphi$ are bijection.
\end{proof}

Let $\bar{\nu};\mu \in Q$ be an $(m| n)$-standard composite partition. We denote by $\Lambda_{\bar{\nu};\mu}$ the corresponding special weight. In case $\mu=0$, we write $\Lambda_{\bar{\nu}}$ for $\Lambda_{\bar{\nu}};0$.
\begin{lem}\label{L:41} Let  $ \Lambda= \Lambda_{\bar{\nu};\mu}\in P_k$ be a special weight corresponding to $\bar{\nu};\mu \in Q_k$ . \\
Set $ \kappa = (\nu_1, \nu_2, \ldots, \nu_k)$ and $\eta = (\nu_{k + 1}, \nu_{k + 2}, \ldots)$. Then we have
\begin{itemize}
\item[(a)] $ \Lambda= (\mu_1, \mu_2, \ldots, \mu_{m - k}, n-\nu_k, \ldots, n - \nu_2, n-\nu_1; -\nu'_n, \ldots, -\nu'_2, -\nu'_1)$;
\item[(b)] $\overline \eta; \mu \in Q_0$ is an $(m-k| n)$-standard composite partition and \\
$\Lambda_{\overline \eta; \mu}= \left (\mu_1, \ldots, \mu_{m-k}; -(\nu_n'-k), -(\nu_{n - 1}'-k), \ldots, -(\nu_1'-k)\right)\in P_0$;
\item [ (c)] $\overline \kappa; 0 \in Q_k$ is an $(k | n)$-standard composite partition and $\Lambda_{\overline \kappa; 0} = (n-\nu_k, \ldots, n -\nu_{1}; -k, \ldots, -k) \in P_k$;
\item[ (d)] $\Lambda = \Lambda_{\overline \eta; \mu} + \Lambda_{\overline \kappa; 0}$;
\end{itemize}
\end{lem}

\begin{proof}[Proof]
 $(a)$ can be readily seen from proof of  Proposition \ref{Pro:2}.\\

Proof of $(b)$. Since $\bar{\nu};\mu \in Q_k$ then 
$$ 
\begin{cases}
    \mu'_{ 1}  \leq m - k  &\mbox{}\\ 
     \eta'_{n} = 0. 
\end{cases} 
 $$
Thus $\overline \eta; \mu \in Q_0$ is an $(m-k| n)$-standard composite partition. It follows from proof of Proposition \ref{Pro:2} that
\begin{align*}
\Lambda_{\overline \eta; \mu}&=  \left(\mu_1, \ldots, \mu_{m-k}; -\eta_n', \ldots, -\eta_1' \right)\\
&=  \left (\mu_1, \ldots, \mu_{m-k}; -(\nu_n'-k)=0, -(\nu_{n - 1}'-k), \ldots, -(\nu_1'-k)\right).
\end{align*}
So $\Lambda_{\overline \eta; \mu} \in P_0$ .\\

We verify $(c)$. We have 
$ \kappa =(\nu_1, \ldots, \nu_k)$ and $\kappa'_n = \nu'_n = k$, thus
$\overline \kappa; 0 \in Q_k$ is an $(k | n)$-standard composite partition. By Proposition \ref{Pro:2}, we have
\begin{align*}
\Lambda_{\overline \kappa; 0}&= (n-\kappa_k, \ldots, n -\kappa_{1}; -\kappa_{n}', \ldots, -\kappa_{1}')\\ 
& = (n-\nu_k, \ldots, n -\nu_{1}; -k, \ldots, -k) \in P_k.
\end{align*}
 
 $(d)$. We have 
\begin{align*}
&\Lambda_{\overline \eta; \mu} + \Lambda_{\overline \kappa; 0}
=(\mu_1\epsilon_1+ \ldots + \mu_{m-k}\epsilon_{m-k} -(\nu_n'-k)\delta_1 -(\nu_{n - 1}'-k)\delta_2 \ldots -(\nu_1'-k)\delta_n)\\
&+ ((n-\nu_k)\epsilon_{m-k+1} \ldots + (n -\nu_{1})\epsilon_m -k \delta_1+ \ldots -k \delta_n)  \\ 
& =\mu_1\epsilon_1+ \ldots + \mu_{m-k}\epsilon_{m-k}+(n-\nu_k)\epsilon_{m-k+1} \ldots + (n -\nu_{1})\epsilon_m
 -(\nu_n')\delta_1 -(\nu_{n - 1}')\delta_2 \ldots -(\nu_1')\delta_n\\
&= \Lambda.
\end{align*}
\end{proof}

\section {Reduction formula for irreducible characters of  $ \frak{gl}(m|n)$}
In this section, we shall establish a reduction formula to represent irreducible characters of $\frak{gl}(m|n)$ in terms of irreducible characters of its subalgebras. Our main ingredient is the character formular of Su-Zhang. We shall start the section by reviewing the formula. We use notations of \cite{Zhang2}
\subsection{Su-Zhang's character formular}
In this subsection we recall a result of  Su-Zhang \cite{Zhang2}, for that we shall need some notations. Assume
\begin{equation}\label{eq:31}
\Lambda = (\lambda_1, \ldots, \lambda_{m_r}, \ldots, \lambda_i, \ldots, \lambda_{m_{1}}, \ldots, \lambda_m; \mu_1, \ldots,\mu_{n_1}, \ldots, \mu_j, \ldots, \mu_{n_r}, \ldots, \mu_n) \in \frak h^* 
\end{equation}
is $r$-fold atypical integral dominant weight with the set of atypical roots
\begin{equation}\label{eq:32}
\Gamma_\Lambda = \{\gamma_1,\ldots, \gamma_r\}.
\end{equation}
Note that we have $ \gamma_1= \epsilon_{m_1}-\delta_{n_1}< \cdots < \gamma_r = \epsilon_{m_r} - \delta_{n_r},$ and $m_r<\cdots<m_1, n_1 <\cdots< n_r$.\\

The total order on $\Delta_1^+$ is defined by
\begin{equation}\label{eq:}
\epsilon_i - \delta_\zeta < \epsilon_j - \delta_\eta  \,\,\mbox{if}\,\,  \zeta - i < \eta - j \,\,\mbox{or}\,\, \zeta - i = \eta - j \,\, \mbox{but} \,\,i > j.
\end{equation}
We call $\gamma_s$ the $s$-th atypical root of $\Lambda$ for $s = 1, 2, \ldots, r$.

 For convenience, we introduce the notation $\Lambda^\rho$ for $\rho$-translation of $\Lambda$ :
\begin{equation}\label{eq:}
\Lambda^\rho = \Lambda+ \rho.
\end{equation}

We define the atypical tuple of $\Lambda$
\begin{equation}\label{eq:33}
aty_\Lambda = (\mu_{n_1}^\rho, \ldots, \mu_{n_r}^\rho) = (\mu_{n_1} - n_1, \ldots, \mu_{n_r} - n_r),
\end{equation}
and we call the $s$-th entry of $aty_\Lambda$ the $s$-th atypical entry of $\Lambda$ for $s = 1, 2, \ldots, r$. We also define the typical tuple of $\Lambda$
\begin{equation}\label{eq:}
typ_\Lambda \in \mathbb{Z}^{m -r| n - r}
\end{equation}
to be the element obtained from $\Lambda^\rho$ by deleting all entries $\lambda_{m_s}^\rho, \mu_{n_s}^\rho$ for $s = 1, 2, \ldots, r$.\\

We call $\Lambda$ lexical if its atypical tuple $ aty_{\Lambda}$ is lexical in the following sence:
$$ \mu_{n_1}^\rho \geq \mu_{n_2}^\rho\geq \ldots \geq \mu_{n_r}^\rho .$$

Corresponding to each atypical root $\gamma_s$  of $\Lambda$, one defines the $\gamma_s$-height of $\Lambda$ by the following formula.
\begin{equation}\label{eq:}
h_s(\Lambda) = \lambda_{m_s} - n_s + s .
\end{equation}

 For $1\leq s \leq t \leq r $, we set
\begin{equation}\label{eq:}
d_{s,t}(\Lambda) = h_t(\Lambda) - h_s(\Lambda) = \lambda_{m_t} - \lambda_{m_s} - n_t + n_s +  t - s.
\end{equation}
 Then $d_{s,t}(\Lambda)$ is non-negative and one can observe that it is the number of integers between the $s$-th atypical entry $\mu_{n_s}^\rho$ and the $\textit{t}$-th atypical entry $\mu_{n_t}^\rho$ which are not entries of $\Lambda^\rho$. In orther words,

\begin{equation}\label{eq:}
d_{s,t}(\Lambda) = \# ([\mu_{n_s}^\rho,\mu_{n_t}^\rho] \texttt{\symbol{92}} Set(\Lambda^\rho)),
\end{equation}
where
\begin{align}
[i,j]&=\left\{ \begin{array}{l}
\{ k \in \mathbb{Z}| i\leq k \leq j\} \quad \mbox{if}\; i \leq j
\\
\emptyset \quad \mbox{otherwise},
\end{array} \right. \mbox{for} \; i,j \in \mathbb{Z}\\
Set(\mu)& = \mbox{the set of the entries of a weight}\,\, \mu.
\end{align}

For $s\leq t$, we say that two atypical roots $\gamma_s, \gamma_t$ of $\Lambda$ are \textit{c-related} if $s = t$ or $d_{s,t}(\Lambda)< t-s$, and are \textit{strongly c-related} if $\gamma_s, \gamma_{u+1}$ are $c$-related for all $r$ such that $s \leq u < t$.\\
The relation $d_{s,t}(\Lambda)< t-s$ is equivalent to\\
\begin{equation}\label{eq:}
\lambda_{m_t} - \lambda_{m_s} < n_t - n_s .
\end{equation} 

A weight $\Lambda$ is said to be totally connected if two atypical roots $\gamma_s, \gamma_t$ of $\Lambda$ are \textit{c-related} for all pairs $(s,t)$ with $s \leq t$.\\

For an $r$-fold atypical weight as given by \eqref{eq:31} and $\sigma$ is an element of symmetric group $S_r$ we define 
$$ \sigma(\Lambda) =  (\lambda_1, \ldots, \lambda_{m_{\sigma^{-1}(r)}}, \ldots, \lambda_i, \ldots, \lambda_{m_{\sigma^{-1}(1)}}, \ldots, \lambda_m; \mu_1, \ldots,\mu_{n_{\sigma^{-1}(1)}}, \ldots, \mu_j, \ldots, \mu_{n_{\sigma^{-1}(r)}}, \ldots, \mu_n).$$
Thus we also have the \textit{dot action}
\begin{equation}\label{eq:34}
\sigma . \Lambda = \sigma(\Lambda + \rho) - \rho.
\end{equation}

Define $S^\Lambda$ to be a subset of the symmetric group $S_r$ consisting of permutations $\sigma$ which do not change the order of $ s < t $ when the atypical roots $\gamma_s$ and $ \gamma_t$ of $\Lambda$ are strongly c-related. That is,
\begin{equation}\label{eq:35}
S^\Lambda = \{\sigma \in  S_r| \sigma^{-1}(s) < \sigma^{-1}(t) \; \mbox{for all }\; s < t \; \mbox {with}\; \gamma_s, \gamma_t \; \mbox{are strongly c-related}\}.
\end{equation}

Let $\Lambda$ in \eqref{eq:31} be an $r$-fold atypical weight with atypical roots ordered as in \eqref{eq:32}: $\gamma_1<\gamma_2<\cdots <\gamma_r$. We define the \textit{normal cone} with vertex $\Lambda$:
\begin{equation}\label{eq:}
C_\Lambda^{Norm} = \{ \Lambda - \sum_{s=1}^r i_s\gamma_s| i_s \geq 0\}.
\end{equation} 

Define partial order $''\leq''$ on $C_\Lambda^{Norm}$ such that for $\lambda, \mu \in C_\Lambda^{Norm}$: $\mu \leq \lambda$ iff $\mu_i \leq \lambda_i$ with $i=1,2,\ldots,m$ and $\mu_{m+j}\geq \lambda_{m+j}$ with $j=1,2,\ldots,n$.

For $\lambda = \Lambda - \sum_{s=1}^r i_s\gamma_s \in C_\Lambda^{Norm} $, we set
$$ |\Lambda - \lambda| = \sum_{s=1}^r i_s\quad.$$
This number is called the level of $\lambda$.
For $\lambda \in C_\Lambda^{Norm}$, denote by $\lambda_{\uparrow}$ the maximal lexical weight which is $\leq \lambda$, namely,
\begin{equation}\label{eq:36}
\lambda _{\uparrow} = max\{ \mu \in C_\Lambda^{Norm}| \mu \leq \lambda, \;\; \mbox{and}\; \; \mu \;\; \mbox{is lexical}\}.
\end{equation} 

Denote by $C_r$ the subset of $S_r$ 
\begin{equation}\label{eq:37}
C_r = \{\pi \in S_r | \pi = (1, 2, \ldots, i_1)(i_1 + 1, i_1 + 2, \ldots, i_1 + i_{-2})\ldots(i_1 + \ldots + i_{t-1}+1, \ldots, r) \},
\end{equation} 
where $i_1, i_2, \ldots, i_t$ are positive integers such that $\sum_{s = 1}^t i_s = r$.
And
\begin{equation}\label{eq:}
\binom{r}{\pi} = \frac{r!}{i_1!i_2!\ldots i_t!}.
\end{equation} 

\begin{thm}\label{T:31}\cite[Theorem 4.9]{Zhang2}
The formal character $ch V(\Lambda)$ of the finite dimensional irreducible $\frak{g}$-module $V(\Lambda)$ is given by
\begin{equation}\label{eq:SZ}
 ch V(\Lambda) = \sum_{\sigma \in S^\Lambda, \pi \in C_r}\frac{1}{r!}\binom{r}{\pi}(-1)^{|\Lambda-(\pi . (\sigma .\Lambda)_{\uparrow})_\uparrow| + l(\pi)} \frac{1}{L_0}\sum_{w \in W}\epsilon (w)w(e^{(\pi . (\sigma .\Lambda)_{\uparrow})_\uparrow + \rho_0}\prod_{\beta \in \Delta_1^+ \backslash{\Gamma_\Lambda}}(1 + e^{-\beta}) ), 
\end{equation}
where $S^\Lambda, \lambda _{\uparrow}, C_r$  are defined above and $L_0 = \prod_{\alpha \in \Delta^+_0}(e^{\alpha/2} - e^{-\alpha/2})$.
\end{thm}

\subsection{Reduction formula for irreducible characters of  $\frak{gl}(m|n)$}
The aim of this subsection is to prove a reduction formula, representing irreducible characters of $\frak{gl}(m|n)$ in terms of irreducible characters of its subalgebras (cf. Theorem \ref{T:ch}).

Denote $\rho(p|q)=(p, p-1,\ldots, 1;-1,-2, \ldots, q)$ with $ p,q$ being positive integers.
\begin{lem}\label{L:51} Let  $ \Lambda=\Lambda_{\bar{\nu};\mu} \in P_k$ be a special weight corresponding to $\bar{\nu};\mu \in Q_k$ . \\
Set $ \kappa = (\nu_1, \nu_2, \ldots, \nu_k)$ and $\eta = (\nu_{k + 1}, \nu_{k + 2}, \ldots)$. Then we have
\begin{itemize}
\item [(a)] $\Gamma_{\Lambda}= \Gamma_{\Lambda_{\overline{\eta}; \mu}}$ and $\Gamma_{\Lambda_{\overline{\kappa}; 0}} = {\emptyset}$;
\item [(b)] $ S^{\Lambda}= S^{\Lambda_{\overline{\eta}; \mu}}$.
\end{itemize}
\end{lem}
\begin{proof}[Proof]
First we prove  $(a)$. From Lemma \ref{L:41} (a), (b) we have 
$$\Lambda= (\mu_1, \mu_2, \ldots, \mu_{m - k}, n-\nu_k, \ldots, n - \nu_2, n-\nu_1; -\nu'_n, \ldots, -\nu'_2, -\nu'_1)$$
 and $$\Lambda_{\overline \eta; \mu}= \left (\mu_1, \ldots, \mu_{m-k}; -(\nu_n'-k), -(\nu_{n - 1}'-k), \ldots, -(\nu_1'-k)\right).$$
Assume $\epsilon_i - \delta_j $, $1 \leq i \leq m-k $ \;\mbox{and}\; $ 1\leq j \leq n$, is a root of $\Lambda$. Then
$$  (\Lambda + \rho(m|n), \epsilon_i - \delta_j) = 0$$
that is
$$  \mu_{i} + m - i + 1 = -(-\nu_{ j}' - j).$$
This is equivalent to
$$  \mu_{i} + (m - k)- i + 1 = -[-(\nu_{ j}' - k) - j].$$
So
$$  (\Lambda_{\overline{\eta}; \mu}+ \rho(m-k|n), \epsilon_i - \delta_j) = 0.$$
Hence, $\epsilon_i - \delta_j $ is a root of $\Lambda_{\overline{\eta}; \mu}$.

On the orther hand, if $m-k <  i \leq m $ and $ 1\leq j \leq n$ then $(\Lambda + \rho(m|n), \epsilon_i - \delta_j) \ne 0$ since $n-\nu_{i} + m - i + 1 < k + 1$ but $-(-\nu_{ j}' - j) \geq k + 1$. Thus, $\epsilon_i - \delta_j $ is not a root of $\Lambda$.  \\
So, we conclude that $$\Gamma_{\Lambda}= \Gamma_{\Lambda_{\overline{\eta}; \mu}}. $$

Next, from property (c) of Lemma \ref{L:41} we have $\Lambda_{\overline \kappa; 0} = (n-\nu_k, \ldots, n -\nu_{1}; -k, \ldots, -k)$. 
Clearly, $\Gamma_{\Lambda_{\overline{\kappa}; 0}} = {\emptyset}$ because of $(n -\nu_1) +1\leq (n -\nu_2) +2\leq \ldots, (n -\nu_k) +k\leq k< -(-k-j)$ for all $j= 1,2,\ldots n$.\\

The second we prove $(b)$. From $(a)$, we have $\Gamma_{\Lambda}= \Gamma_{\Lambda_{\overline{\eta}; \mu}}$. By assumption, $ 1 \leq m_t \leq m_s \leq m - k$, it is easy to see that $  \gamma_s, \gamma_t \in \Gamma_{\Lambda}$ are $c$- related if only if $\gamma_s, \gamma_t \in \Gamma_{\Lambda_{\overline{\eta}; \mu}}$ are $c$ - related. Thus $ S^{\Lambda}= S^{\Lambda_{\overline{\eta}; \mu}}$. 
\end{proof}

\begin{lem}\label{L:52} Let  $ \Lambda \in P_k$ be a special weight corresponding to $\bar{\nu};\mu \in Q_k$ . \\
Set $ \kappa = (\nu_1, \nu_2, \ldots, \nu_k)$ and $\eta = (\nu_{k + 1}, \nu_{k + 2}, \ldots)$. Then we have
$$  (\pi . (\sigma . \Lambda)_{\uparrow})_{\uparrow}  =  (\pi . (\sigma . \Lambda_{\overline{\eta}; \mu})_{\uparrow})_{\uparrow} + \Lambda_{\overline \kappa; 0}$$ and
 $$|\Lambda - (\pi . (\sigma . \Lambda)_{\uparrow})_{\uparrow} | = |\Lambda_{\overline{\eta}; \mu}- (\pi . (\sigma . \Lambda_{\overline{\eta}; \mu})_{\uparrow})_{\uparrow}|,$$
where $\pi \in C_r , \sigma \in S^{\Lambda}$.
\end{lem}
\begin{proof}[Proof]
This follows from Lemma \ref{L:41} and Lemma \ref{L:51} with notice that $\sigma . \Lambda= \sigma . \Lambda_{\overline{\eta}; \mu} + \Lambda_{\overline \kappa; 0}$ and $\pi . (\sigma . \Lambda)_{\uparrow}= \pi . (\sigma . \Lambda_{\overline{\eta}; \mu})_{\uparrow} + \Lambda_{\overline \kappa; 0}$.
\end{proof}
\begin{lem}\label{L:rho} Let $m, n, p, q$ be nonnegative integers such that $m = p + q$. We have
\begin{itemize}
\item [(a)] \begin{align*}
L_0(p|n) &\prod_{0 \leq i < j \leq q}(e^{(\epsilon_{p + i}- \epsilon_{p + j})/2}- e^{-(\epsilon_{p + i}- \epsilon_{p + j})/2}) \prod_{i= 1}^p \prod_{j = 1}^q (e^{\epsilon_{i}}- e^{\epsilon_{p+ j}})\\ 
& = L_0(m|n) e^{\frac{1}{2}(q, \ldots, q, p, \ldots, p; 0, \ldots, 0)},
\end{align*}
where $$L_0(m|n) = \prod_{\alpha \in \Delta^+_0}(e^{\alpha/2} - e^{-\alpha/2})$$ and
$$L_0(p|n) = \prod_{\alpha \in \Delta^+_0\setminus \{(\epsilon_i -\epsilon_j) | p + 1 \leq  j \leq m\}}(e^{\alpha/2} - e^{-\alpha/2}).$$
\item [(b)] 
\begin{align*}
\rho_0(m|n) = &\rho_0(p|n) + \frac{1}{2}(0, \ldots, 0, q -1, q - 3, \ldots, 1 - q; 0, \ldots, 0)\\ 
& +  \frac{1}{2}(q, \ldots, q, -p, \ldots, -p; 0, \ldots, 0),
\end{align*}
where $$\rho_0(m| n) = \frac{1}{2}(m-1, m-3, \ldots, 1-m; n - 1, n - 3, \ldots, 1 - n)$$ and 
$$\rho_0(p| n) = \frac{1}{2}((m - q)-1, (m - q)-3, \ldots, 1-(m -q); n - 1, n - 3, \ldots, 1 - n).$$
\end{itemize}
\end{lem}

\begin{proof}
 We have
\begin{align*}
\prod_{i= 1}^p \prod_{j = 1}^q (e^{\epsilon_{i}}- e^{\epsilon_{p+ j}})&= \prod_{i= 1}^p \prod_{j = 1}^q e^{\epsilon_i/2}e^{\epsilon_{p+j}/2}(e^{(\epsilon_{i}- \epsilon_{p + j})/2}- e^{-(\epsilon_{ i}- \epsilon_{p + j})/2})\\ 
& = \left(\prod_{i= 1}^p e^{\epsilon_i/2}\right)^q  \left(\prod_{j = 1}^q e^{\epsilon_{p+j}/2}\right)^p \prod_{i= 1}^p \prod_{j = 1}^q (e^{(\epsilon_{i}- \epsilon_{p + j})/2}- e^{-(\epsilon_{ i}- \epsilon_{p + j})/2}).
\end{align*}
Thus
\begin{align*}
&L_0(p|n) \prod_{0 \leq i < j \leq q}(e^{(\epsilon_{p + i}- \epsilon_{p + j})/2}- e^{-(\epsilon_{p + i}- \epsilon_{p + j})/2}) \prod_{i= 1}^p \prod_{j = 1}^q (e^{\epsilon_{i}}- e^{\epsilon_{p+ j}})\\ 
&=  L_0(m|n) \left(\prod_{i= 1}^p e^{\epsilon_i/2}\right)^q  \left(\prod_{j = 1}^q e^{\epsilon_{p+j}/2}\right)^p \\
& = L_0(m|n) e^{\frac{1}{2}(q, \ldots, q, p, \ldots, p; 0, \ldots, 0)}.
\end{align*}
This show (a).\\

Finally, (b) is obvious.
\end{proof}

\begin{lem}\label{L:kapa}
Let $\Lambda = (\lambda_1, \ldots, \lambda_m; q,\ldots, q)$  be an integral dominant weight  such that $\Gamma_{\Lambda} = {\emptyset}$. Then
$$ 
ch\; V(\Lambda) = \frac{\prod_{\beta\in\Delta_1^+} (1 + e^{-\beta})}{\prod_{i < j}(e^{(\epsilon_{ i}- \epsilon_{ j})/2}- e^{-(\epsilon_{ i}- \epsilon_{ j})/2})} \sum_{w \in S_m}\epsilon(w) w(e^{\Lambda + \frac{1}{2}(m-1, m-3, \ldots, 1-m; 0,\ldots, 0)}).
 $$
\end{lem}

\begin{proof} 
Apply Theorem \ref{T:31} with $r=0$ we have
\begin{align}
ch \;V(\Lambda) 
& =\frac{1}{L_0}
 \sum_{w \in S_m\times S_n} \epsilon(w) w\left(e^{ \Lambda + \rho_0 } \prod_{\beta\in\Delta_1^+ }(1 + e^{-\beta})\right) \\ 
 &=\frac{\prod_{\beta\in\Delta_1^+ }(1 + e^{-\beta})}{L_0}
 \sum_{w \in S_m\times S_n} \epsilon(w) w(e^{ \Lambda + \rho_0 })\\ 
& =  \frac{\prod_{\beta\in\Delta_1^+}(1 + e^{-\beta})}{\prod_{0\leq i < j \leq m}(e^{(\epsilon_{ i}- \epsilon_{ j})/2}- e^{-(\epsilon_{ i}- \epsilon_{ j})/2})}
\times \frac{\sum_{w \in S_m\times S_n} \epsilon(w) w(e^{ \Lambda + \rho_0 })}{\prod_{1\leq i < j \leq n}(e^{(\delta_{ i}- \delta_{ j})/2}- e^{-(\delta_{ i}- \delta_{ j})/2})}\\
& =  \frac{\prod_{\beta\in\Delta_1^+ }(1 + e^{-\beta})}{\prod_{0\leq i < j \leq m}(e^{(\epsilon_{ i}- \epsilon_{ j})/2}- e^{-(\epsilon_{ i}- \epsilon_{ j})/2})}
\times \frac{\sum_{w_0 \in S_m} \epsilon(w_0) w_0(\sum_{w_1 \in S_n} \epsilon(w_1) w_1(e^{ \Lambda + \rho_0 }))}{\prod_{1\leq i < j \leq n}(e^{(\delta_{ i}- \delta_{ j})/2}- e^{-(\delta_{ i}- \delta_{ j})/2})},
\end{align}
where $w= w_0\times w_1 \in S_m\times S_n$.

Since $\Lambda = (\lambda_1, \ldots, \lambda_m; q,\ldots, q)$ then $ w_1(e^{ \Lambda +  \frac{1}{2}(m-1, m-3, \ldots, 1-m; 0,\ldots, 0) })=  e^{\Lambda +  \frac{1}{2}(m-1, m-3, \ldots, 1-m; 0,\ldots, 0) } $, for all $w_1 \in S_n$. And $ \rho_0=\frac{1}{2}(m-1, m-3, \ldots, 1-m; 0,\ldots, 0) +  \frac{1}{2}(0,\ldots, 0;n - 1, n- 3, \ldots, 1-n) $. These follow

\begin{align}
ch \;V(\Lambda) 
&=  \frac{\prod_{\beta\in\Delta_1^+ }(1 + e^{-\beta})}{\prod_{0\leq i < j \leq m}(e^{(\epsilon_{ i}- \epsilon_{ j})/2}- e^{-(\epsilon_{ i}- \epsilon_{ j})/2})}\\
&\times \frac{\sum_{w_0 \in S_m} \epsilon(w_0) w_0(e^{ \Lambda +  \frac{1}{2}(m-1, m-3, \ldots, 1-m; 0,\ldots, 0) }\sum_{w_1 \in  S_n} \epsilon(w_1) w_1(e^{ \frac{1}{2}(0,\ldots, 0;n - 1, n- 3, \ldots, 1-n) }))}{\prod_{1\leq i < j \leq n}(e^{(\delta_{ i}- \delta_{ j})/2}- e^{-(\delta_{ i}- \delta_{ j})/2})}.
\end{align}  

Because $ \prod_{1\leq i < j \leq n}(e^{(\delta_{ i}- \delta_{ j})/2}- e^{-(\delta_{ i}- \delta_{ j})/2})= \prod_{1\leq i < j \leq n}(e^{\delta_{ i}}- e^{ \delta_{ j}})(\prod_{i = 1}^n e^{\delta_i})^{-(n-1)/2} $, so
\begin{align}
ch \;V(\Lambda) 
&=  \frac{\prod_{\beta\in\Delta_1^+ }(1 + e^{-\beta})}{\prod_{0\leq i < j \leq m}(e^{(\epsilon_{ i}- \epsilon_{ j})/2}- e^{-(\epsilon_{ i}- \epsilon_{ j})/2})}\\
&\times \frac{\sum_{w_0 \in S_m} \epsilon(w_0) w_0(e^{ \Lambda +  \frac{1}{2}(m-1, m-3, \ldots, 1-m; 0,\ldots, 0) }\sum_{w_1 \in  S_n} \epsilon(w_1) w_1(e^{ (0,\ldots, 0;n - 1, n- 2, \ldots, 0) }))}{\prod_{1\leq i < j \leq n}(e^{\delta_{ i}}- e^{ \delta_{ j}})}.
\end{align}  

Since 
$$ \sum_{w_1 \in  S_n} \epsilon(w_1) w_1(e^{ (0,\ldots, 0;n - 1, n- 2, \ldots, 0) })) = \prod_{1\leq i < j \leq n}(e^{\delta_{ i}}- e^{ \delta_{ j}}), $$
then
\begin{align*}
ch \;V(\Lambda) =  \frac{\prod_{\beta\in\Delta_1^+ }(1 + e^{-\beta})}{\prod_{0\leq i < j \leq m}(e^{(\epsilon_{ i}- \epsilon_{ j})/2}- e^{-(\epsilon_{ i}- \epsilon_{ j})/2})} \times \sum_{w_0 \in S_m} \epsilon(w_0) w_0(e^{ \Lambda +  \frac{1}{2}(m-1, m-3, \ldots, 1-m; 0,\ldots, 0) }).
\end{align*}  
\end{proof}

We first define some substitution rules.
Let $x=(x_1, x_2, \ldots,x_m)$ and $ y=(y_1, y_2, \ldots, y_n)$ be the sets of variables. For an any weight $\lambda = (\alpha_1, \alpha_2, \ldots, \alpha_m; \beta_1, \beta_2, \ldots, \beta_n)$, we denote by $(x;y)^{\lambda}$ the monomial $x_1^{\alpha_1}\ldots x_2^{\alpha_m}y_1^{\beta_1} \ldots y_n^{\beta_n}$.\\
Let $\underline{r} := \{r_1, r_2, \ldots, r_{m-k}\} \subset \{1, 2, \ldots, m\}$ and $ \underline{s}:= \{r_{m-k + 1},r_{m-k + 2}, \ldots, r_{m}\} = \{1, 2, \ldots, m\} \backslash \underline{r}$, where $k$ is a nonnegative integer less than or equal $m$.
For $f $, a rational function in $x_1, x_2, \ldots,x_{m-k}, y_1, y_2, \ldots, y_n$, define $ \chi_{\underline{r}}(f) $ to be the rational function obtained from $f$ by substituting $x_i$ by $x_{r_{i}}$. That is, if
$$ f = \frac{P(x_1, x_2, \ldots,x_{m-k}, y_1, y_2, \ldots, y_n)}{Q(x_1, x_2, \ldots,x_{m-k}, y_1, y_2, \ldots, y_n)} $$
then
$$  \chi_{\underline{r}}(f) = \frac{P(x_{r_{1}},\ldots, x_{r_{m-k}}, y_1, \ldots, y_n)}{Q(x_{r_{1}},\ldots, x_{r_{m-k}}, y_1, \ldots, y_n)}.$$

\begin{thm}\label{T:ch}
Let  $ \Lambda $ be a special weight in $P_k$ and let $\bar{\nu};\mu \in Q_k$ be the corresponding composite partition. Then
\begin{equation}\label{eq:ch}
{\rm ch} V\; (\Lambda) = \sum_{\underline{r}, \underline{s}}
\frac{(\prod_{i = 1}^{m-k}  x_{r_i})^k  \chi_{\underline{r}}(ch\;V(\Lambda_{\overline\eta;\mu}) ) \chi_{\underline{s}}(ch\;V(\Lambda_{\overline\kappa}))}{\prod_{i= 1}^{m-k} \prod_{j = 1}^k (x_{r_i}- x_{r_{m-k+j}})},
\end{equation}
where $\nu=  (\nu_1, \nu_2, \ldots, \nu_k,\nu_{k + 1}, \nu_{k + 2}, \ldots), \kappa = (\nu_1, \nu_2, \ldots, \nu_k)$, $\eta = (\nu_{k + 1}, \nu_{k + 2}, \ldots)$ and the sum is over all possible decomposition $\{1, 2, \ldots, m\} =\underline{r} \cup \underline{s}$ with $|\underline{r}|= m-k$, $|\underline{s}|= k$.
\end{thm}

\begin{proof}[Proof ]
Because $ \kappa = (\nu_1, \nu_2, \ldots, \nu_k)$ so we have
$\Lambda_{\overline\kappa}=(n-\nu_k,n-\nu_{k-1},\ldots, n-\nu_1; -k,\ldots,-k)$ (because of property (c) of Lemma \ref{L:41})
and $\Gamma_{\Lambda_{\kappa}}=\emptyset$ (because of property (a) of Lemma \ref{L:51}).
From Lemma \ref{L:kapa}, we have
\begin{align}\label{Eq:55}
ch V(\Lambda_{\overline\kappa}) &= \frac{\prod_{\beta\in\Delta_1^+(k|n) }(1 + e^{-\beta})}{\prod_{m-k+1\leq i < j \leq m}(e^{(\epsilon_{ m-k+i}- \epsilon_{  m-k+j})/2}- e^{-(\epsilon_{ m-k+i}- \epsilon_{ m-k+j})/2})}\\
&\times \sum_{w_0 \in S_k}\epsilon(w_0) w_0(e^{\Lambda_{\overline\kappa} + \frac{1}{2}((k-1)\epsilon_{m-k +1}+ (k-3)\epsilon_{m-k +2}+ \ldots+ (1-k)\epsilon_{m})}).
 \end{align}

We set
\begin{align}
f: &= \frac{(\prod_{i = 1}^{m-k}  e^{\epsilon_i})^k  ch\;V(\Lambda_{\overline\eta;\mu}) ch\;V(\Lambda_{\overline\kappa})}{\prod_{i= 1}^{m-k} \prod_{j = 1}^k (e^{\epsilon_i}- e^{\epsilon_{m-k +j}})}\\ 
& =  \frac{(\prod_{i = 1}^{m-k}  e^{\epsilon_i})^k  }{\prod_{i= 1}^{m-k} \prod_{j = 1}^k (e^{\epsilon_i}- e^{\epsilon_{m-k +j}})}
\times ch\;V(\Lambda_{\overline\eta;\mu})  \times ch\;V(\Lambda_{\overline\kappa}).
\end{align}
By using Theorem \ref{T:31} and Eq. \eqref{Eq:55}, we obtain
\begin{align}
f&=  \frac{(\prod_{i = 1}^{m-k}  e^{\epsilon_i})^k  }{\prod_{i= 1}^{m-k} \prod_{j = 1}^k (e^{\epsilon_i}- e^{\epsilon_{m-k +j}})} \\
 &\times \sum _{\sigma \in S^{\Lambda_{\overline{\eta}; \mu}}, \pi \in C_r}\frac{1}{r!}\binom{r}{\pi}(-1)^{|\Lambda_{\overline{\eta}; \mu}- (\pi . (\sigma . \Lambda_{\overline{\eta}; \mu})_{\uparrow})_{\uparrow}| + l(\pi)}\\
&\times \frac{1}{L_0(m-k|n)}
 \sum_{w \in S_{m-k}\times S_n} \epsilon(w) w\left(e^{(\pi . (\sigma . \Lambda_{\overline{\eta}; \mu})_{\uparrow})_{\uparrow} + \rho_0(m-k|n) } \prod_{\beta\in\Delta_1^+(m-k|n) \backslash{\Gamma_{\Lambda_{\overline{\eta}; \mu}}}}(1 + e^{-\beta})\right) \\
& \times  \frac{\prod_{\beta\in\Delta_1^+(k|n) }(1 + e^{-\beta})}{\prod_{m-k+1\leq i < j \leq m}(e^{(\epsilon_{ m-k+i}- \epsilon_{  m-k+j})/2}- e^{-(\epsilon_{ m-k+i}- \epsilon_{ m-k+j})/2})} \\
& \times \sum_{w_0 \in S_k}\epsilon(w_0) w_0(e^{\Lambda_{\overline\kappa} + \frac{1}{2}((k-1)\epsilon_{m-k +1}+ (k-3)\epsilon_{m-k +2}+ \ldots+ (1-k)\epsilon_{m})}),
\end{align}
where $\Delta_1^+(m-k|n)= \{\epsilon_i - \delta_j | 1 \leq i \leq m-k, 1 \leq j \leq n \}  $ and $\Delta_1^+(k|n)=\{\epsilon_{m-k+i} - \delta_j | 1 \leq i \leq k, 1 \leq j \leq n \} $.

  From property $(a)$ of Lemma \ref{L:rho}, we deduce
\begin{align}
f&=  \frac{(\prod_{i = 1}^{m-k}  e^{\epsilon_i})^k  }{ L_0(m|n) e^{\frac{1}{2}(k, \ldots, k, m-k, \ldots, m-k; 0, \ldots, 0)}}\\
 &\times \sum _{\sigma \in S^{\Lambda_{\overline{\eta}; \mu}}, \pi \in C_r}\frac{1}{r!}\binom{r}{\pi}(-1)^{|\Lambda_{\overline{\eta}; \mu}- (\pi . (\sigma . \Lambda_{\overline{\eta}; \mu})_{\uparrow})_{\uparrow}| + l(\pi)}\\
&\times \sum_{w \in S_{m-k}\times S_n} \epsilon(w) w(e^{(\pi . (\sigma . \Lambda_{\overline{\eta}; \mu})_{\uparrow})_{\uparrow} + \rho_0(m-k|n) } \prod_{\beta\in\Delta_1^+(m-k|n) \backslash{\Gamma_{\Lambda_{\overline{\eta}; \mu}}}}(1 + e^{-\beta}))\\
& \times \prod_{\beta\in\Delta_1^+(k|n)}(1 + e^{-\beta}) \sum_{w_0 \in S_k}\epsilon(w_0) w_0(e^{\Lambda_{\overline\kappa} + \frac{1}{2}((k-1)\epsilon_{m-k +1}+ (k-3)\epsilon_{m-k +2}+ \ldots + (1-k)\epsilon_{m})}).
\end{align}
We can rewrite
\begin{align}
f&=   \frac{e^{\frac{1}{2}(k, \ldots, k, -(m-k), \ldots,- (m-k); 0, \ldots, 0)} }{ L_0(m|n) }\\
&\times \sum _{\sigma \in S^{\Lambda_{\overline{\eta}; \mu}}, \pi \in C_r}\frac{1}{r!}\binom{r}{\pi}(-1)^{|\Lambda_{\overline{\eta}; \mu}- (\pi . (\sigma . \Lambda_{\overline{\eta}; \mu})_{\uparrow})_{\uparrow}| + l(\pi)}\\
&\times \sum_{w \in S_{m-k}\times S_n} \epsilon(w) w\left(e^{(\pi . (\sigma . \Lambda_{\overline{\eta}; \mu})_{\uparrow})_{\uparrow} + \rho_0(m-k|n) } \prod_{\beta\in\Delta_1^+(m|n) \backslash{\Gamma_{\Lambda_{\overline{\eta}; \mu}}}}(1 + e^{-\beta})\right)\\
& \times \sum_{w_0 \in S_k}\epsilon(w_0) w_0(e^{\Lambda_{\overline\kappa} + \frac{1}{2}((k-1)\epsilon_{m-k +1}+ (k-3)\epsilon_{m-k +2}+ \ldots+ (1-k)\epsilon_{m})}),
\end{align}
where $\Delta_1^+(m|n)=\Delta_1^+$.

We imply
\begin{align}
f&=  \frac{e^{\frac{1}{2}(k, \ldots, k, -(m-k), \ldots,-(m-k); 0, \ldots, 0)} }{ L_0(m|n) }\\
& \times \sum _{\sigma \in S^{\Lambda_{\overline{\eta}; \mu}}, \pi \in C_r}\frac{1}{r!}\binom{r}{\pi}(-1)^{|\Lambda_{\overline{\eta}; \mu}- (\pi . (\sigma . \Lambda_{\overline{\eta}; \mu})_{\uparrow})_{\uparrow}| + l(\pi)}\\
& \times \sum_{w \in (S_{m-k}\times S_k)\times S_n} \epsilon(w) w \left( e^{(\pi . (\sigma . \Lambda_{\overline{\eta}; \mu})_{\uparrow})_{\uparrow} + \rho_0(m-k|n) } e^{\Lambda_{\overline\kappa} + \frac{1}{2}((k-1)\epsilon_{m-k +1}+ (k-3)\epsilon_{m-k +2}+ \ldots+ (1-k)\epsilon_{m})} \right.\\
&\left. \times \prod_{\beta\in\Delta_1^+(m|n) \backslash{\Gamma_{\Lambda_{\overline{\eta}; \mu}}}}(1 + e^{-\beta}) \right).
\end{align}

By using Lemma \ref{L:51} and Lemma \ref{L:52} we have
\begin{align}
f&=    \frac{e^{\frac{1}{2}(k, \ldots, k, -(m-k), \ldots,- (m-k); 0, \ldots, 0)} }{ L_0(m|n) }\\
&\sum _{\sigma \in S^{\Lambda}, \pi \in C_r}\frac{1}{r!}\binom{r}{\pi}(-1)^{|\Lambda- (\pi . (\sigma . \Lambda)_{\uparrow})_{\uparrow}| + l(\pi)}\\
&\times \sum_{w \in (S_{m-k}\times S_k)\times S_n} \epsilon(w) w\left(e^{(\pi . (\sigma . \Lambda)_{\uparrow})_{\uparrow} + \rho_0(m-k|n) + \frac{1}{2}((k-1)\epsilon_{m-k +1}+ (k-3)\epsilon_{m-k +2}+ \ldots + (1-k)\epsilon_{m}) } \right.\\
& \left. \times \prod_{\beta\in\Delta_1^+(m|n) \backslash{\Gamma_{\Lambda}}}(1 + e^{-\beta}) \right).
\end{align}

It follows from property $(b)$ of Lemma \ref{L:rho}:
\begin{align}\label{Eq:f}
f&=  \sum _{\sigma \in S^{\Lambda}, \pi \in C_r}\frac{1}{r!}\binom{r}{\pi}(-1)^{|\Lambda- (\pi . (\sigma . \Lambda)_{\uparrow})_{\uparrow}| + l(\pi)}\\
& \frac{1 }{ L_0(m|n) } \sum_{w \in (S_{m-k}\times S_k)\times S_n} \epsilon(w) w\left(e^{(\pi . (\sigma . \Lambda)_{\uparrow})_{\uparrow} + \rho_0(m|n)  }  \prod_{\beta\in\Delta_1^+(m|n) \backslash{\Gamma_{\Lambda}}}(1 + e^{-\beta})\right).
\end{align}
Now we consider the right of \eqref{eq:ch}. Set
 \begin{align}
 Rhs &:=  \sum_{\underline{r}, \underline{s}}\frac{(\prod_{i = 1}^{m-k}  x_{r_i})^k  \chi_{\underline{r}}(ch\;V(\Lambda_{\overline\eta;\mu}) ) \chi_{\underline{s}}(ch\;V(\Lambda_{\overline\kappa}))}{\prod_{i= 1}^{m-k} \prod_{j = 1}^k (x_{r_i}- x_{r_{m-k+j}})}\\
& = \sum_{\underline{r}, \underline{s}}\chi_{\underline{r}\cup \underline{s}}(f).
\end{align}

By applying substitution rule to \eqref{Eq:f} we have
 \begin{align}
 Rhs &:= \sum_{\underline{r}, \underline{s}} \chi_{\underline{r}\cup \underline{s}}(f)\\
&=\sum_{\underline{r}, \underline{s}}\left(
 \sum _{\sigma \in S^{\Lambda}, \pi \in C_r}\frac{1}{r!}\binom{r}{\pi}(-1)^{|\Lambda- (\pi . (\sigma . \Lambda)_{\uparrow})_{\uparrow}| + l(\pi)} \frac{(\prod_{i = 1}^{m}  x_{r_i})^{\frac{m-1}{2}}(\prod_{j= 1}^{n}  y_{j})^{\frac{n-1}{2}} }{ \prod_{1\leq i < j\leq m} ( x_{r_i}- x_{r_j} )\prod_{1\leq i < j\leq n} ( y_{i}- y_{j} )}\right.\\
&\left.  \sum_{w \in (S_{m-k}\times S_k)\times S_n} \epsilon(w) w \left((x';y)^{(\pi . (\sigma . \Lambda)_{\uparrow})_{\uparrow} + \rho_0(m|n)  }  \prod_{\beta\in\Delta_1^+(m|n) \backslash{\Gamma_{\Lambda}}}(1 +(x';y)^{-\beta})\right)  \right),
\end{align}
where $x' = (x_{r_1},x_{r_2},\ldots, x_{r_m})$.
Hence

\begin{align}
 Rhs &:=\sum _{\sigma \in S^{\Lambda}, \pi \in C_r}\frac{1}{r!}\binom{r}{\pi}(-1)^{|\Lambda- (\pi . (\sigma . \Lambda)_{\uparrow})_{\uparrow}| + l(\pi)} \sum_{\underline{r}, \underline{s}}\left(
  \frac{(\prod_{i = 1}^{m}  x_{i})^{\frac{m-1}{2}}(\prod_{j= 1}^{n}  y_{j})^{\frac{n-1}{2}} }{ \prod_{1\leq i < j\leq m} ( x_{r_i}- x_{r_j} )\prod_{1\leq i < j\leq n} ( y_{i}- y_{j} )}\right.\\
&\left.  \sum_{w \in (S_{m-k}\times S_k)\times S_n} \epsilon(w) w\left((x';y)^{(\pi . (\sigma . \Lambda)_{\uparrow})_{\uparrow} + \rho_0(m|n)  }  \prod_{\beta\in\Delta_1^+(m|n) \backslash{\Gamma_{\Lambda}}}(1 +(x';y)^{-\beta})\right)  \right).
\end{align}

For a fixed pair $\underline{r}, \underline{s}$, we have

\begin{align}
 &
  \frac{(\prod_{i = 1}^{m}  x_{i})^{\frac{m-1}{2}}(\prod_{j= 1}^{n}  y_{j})^{\frac{n-1}{2}} }{ \prod_{1\leq i < j\leq m} ( x_{r_i}- x_{r_j} )\prod_{1\leq i < j\leq n} ( y_{i}- y_{j} )}\\
&\times \sum_{w \in (S_{m-k}\times S_k)\times S_n} \epsilon(w) w\left((x';y)^{(\pi . (\sigma . \Lambda)_{\uparrow})_{\uparrow} + \rho_0(m|n)  }  \prod_{\beta\in\Delta_1^+(m|n) \backslash{\Gamma_{\Lambda}}}(1 +(x';y)^{-\beta})\right)
\end{align}
is equal to
\begin{align}
 &  \frac{(\prod_{i = 1}^{m}  x_{i})^{\frac{m-1}{2}}(\prod_{j= 1}^{n}  y_{j})^{\frac{n-1}{2}} }{ \prod_{1\leq i < j\leq m} ( x_{i}- x_{j} )\prod_{1\leq i < j\leq n} ( y_{i}- y_{j} )}\epsilon(\tau)\\
& \times \sum_{i =1}^{(m-k)!(k)!}\sum_{w_1 \in  S_n} \epsilon(w^i\times w_1) w^i\times w_1\left((x;y)^{(\pi . (\sigma . \Lambda)_{\uparrow})_{\uparrow} + \rho_0(m|n)  }  \prod_{\beta\in\Delta_1^+(m|n) \backslash{\Gamma_{\Lambda}}}(1 +(x;y)^{-\beta})\right),
\end{align}
where $w^1, w^2, \ldots, w^{k!(m-k)!}$ are $k!(m-k)!$ the different permutations of  $\{1,2,\ldots, m\}$. We need a little explanation for the $w^i$. If $\{1, 2, \ldots, m-k\}\backslash \underline{r} = \{u_1, u_2,\ldots, u_h\}$ then $\{m-k+1, \ldots, m\}\backslash \underline{s}=\{v_1, v_2, \ldots, v_h\}$. We put $\tau = (u_1v_1)(u_2v_2) \ldots (u_hv_h)$ is the permutation of $\{1,2,\ldots, m\}$. This permutation is a product of cycles which length equal to 2. Assume that $\sigma_1$ (resp. $\sigma_2$) is a permutation of $\{r_1, r_2, \ldots, r_{m-k}\}$ (resp. $\{r_{m-k + 1},r_{m-k + 2}, \ldots, r_{m}\}$) then $w^i, i=1, 2, \ldots, k!(m-k)!$ will be of the form $\sigma_1 \sigma_2\tau $ (cf. Example \ref{Ex}). It's easy to see that
$$  \frac{1}{ \prod_{1\leq i < j\leq m} ( x_{r_i}- x_{r_j} )} =  \frac{1}{ \prod_{1\leq i < j\leq m} ( x_{i}- x_{j} )}\epsilon(\tau).$$

Because $\underline{r}$ is a set containing $k$ different elements of $\{1,2,\ldots, m\}$ so there are $\frac{m!}{k!(m-k!)}$ different sets $\underline{r}$ (in the orther words, there are $\frac{m!}{k!(m-k!)}$ different pairs of  $\underline{r}, \underline{s}$). Thus 
\begin{align}
 Rhs = &\sum _{\sigma \in S^{\Lambda}, \pi \in C_r}\frac{1}{r!}\binom{r}{\pi}(-1)^{|\Lambda- (\pi . (\sigma . \Lambda)_{\uparrow})_{\uparrow}| + l(\pi)}\frac{(\prod_{i = 1}^{m}  x_{i})^{\frac{m-1}{2}}(\prod_{j= 1}^{n}  y_{j})^{\frac{n-1}{2}} }{ \prod_{1\leq i < j\leq m} ( x_{i}- x_{j} )\prod_{1\leq i < j\leq n} ( y_{i}- y_{j} )}\\
& \sum_{w \in S_{m}\times S_n} \epsilon(w) w\left((x;y)^{(\pi . (\sigma . \Lambda)_{\uparrow})_{\uparrow} + \rho_0(m|n)  }  \prod_{\beta\in\Delta_1^+(m|n) \backslash{\Gamma_{\Lambda}}}(1 +(x;y)^{-\beta})\right)\\
&= ch\; V(\Lambda),
\end{align}
where the last equality comes from Theorem \ref{T:31}:
\begin{align}
ch\; V(\Lambda) &=    \sum _{\sigma \in S^{\Lambda}, \pi \in C_r}\frac{1}{r!}\binom{r}{\pi}(-1)^{|\Lambda- (\pi . (\sigma . \Lambda)_{\uparrow})_{\uparrow}| + l(\pi)} \frac{1 }{ L_0(m|n) }\\
&\times   \sum_{w \in S_m\times S_n} \epsilon(w) w(e^{(\pi . (\sigma . \Lambda)_{\uparrow})_{\uparrow} + \rho_0(m|n)  }  \prod_{\beta\in\Delta_1^+(m|n) \backslash{\Gamma_{\Lambda}}}(1 + e^{-\beta})).
\end{align}

\end{proof}

\begin{Exa} \label{Ex} Assume that $\Lambda = (1,0,-1;-1,-2) \in P_k $ is special weight of $\frak{gl(3|2)}$. Then
\begin{equation}\label{Exa}
{\rm ch} V\; (\Lambda)  = \sum_{\underline{r}}\frac{ x_{r_1} x_{r_2} \chi_{\underline{r}}(ch\;V(\Lambda_{\overline\eta;\mu}) ) \chi_{\underline{s}}(ch\;V(\Lambda_{\overline\kappa}))}{ (x_{r_1}- x_{r_{3}}) (x_{r_2}- x_{r_{3}})},
\end{equation}
where $\Lambda_{\overline\eta;\mu}=(1,0;0,-1), \Lambda_{\overline\kappa}=(-1;-1,-1)$, $\underline{r} = \{r_1, r_2, \} \subset \{1, 2, 3\}$ and $ \underline{s}= \{r_{3} \}= \{1, 2, 3\}\backslash \underline{r} $.
\end{Exa}

\begin{proof}
From Lemma \ref{L:51}, we have
\begin{itemize}
\item [(a)] $\Gamma_{\Lambda}= \Gamma_{\Lambda_{\overline{\eta}; \mu}}=\{\gamma_1 = \epsilon_1- \delta_2, \gamma_2= \epsilon_2 - \delta_1 \}$ ;
\item [ (b)] $ S^{\Lambda}= S^{\Lambda_{\overline{\eta}}}= \emptyset$ because $\gamma_1, \gamma_2 \in \Gamma_{\Lambda}$ are not $c$- related.
\end{itemize}
 Because $r=2$ so $ C_2 = S_2 $ is the symmetric group.

Since $\pi \in S_2 $ then $\pi = (1)$ or $\pi = (12)$.\\
If $\pi = (1)$ then $(\pi . \Lambda_{\overline{\eta}; \mu})_{\uparrow}=\Lambda_{\overline{\eta}; \mu}$ and $(\pi . \Lambda)_{\uparrow}=\Lambda$.\\
If $\pi = (12)$ then  $(\pi . \Lambda_{\overline{\eta}; \mu})_{\uparrow}=(-1,0;0,1)$ and $(\pi . \Lambda)_{\uparrow}=(-1,0, -1;-1,0)$.

We need to compute $ \chi_{\underline{r}}(ch\;V(\Lambda_{\overline\eta;\mu}) ), \chi_{\underline{s}}(ch\;V(\Lambda_{\overline\kappa}))$ and $\chi_{\underline{t}}({\rm ch} V\; (\Lambda) )$.

First, we compute $\chi_{\underline{r}}(ch\;V(\Lambda_{\overline\eta;\mu}) )$. We have
\begin{align*}
 {\rm ch} V\; (\Lambda_{\overline\eta; \mu}) &= \sum _{ \pi \in S_2}\frac{1}{2!}\binom{2}{\pi}(-1)^{|\Lambda_{\overline{\eta}; \mu}- (\pi . \Lambda_{\overline{\eta}; \mu})_{\uparrow}| + l(\pi)}\\
&\times \frac{1}{L_0(2|2)}
 \sum_{w \in S_{2}\times S_2} \epsilon(w) w\left(e^{(\pi .  \Lambda_{\overline{\eta}; \mu})_{\uparrow} + \rho_0(2|2) } \prod_{\beta\in\Delta_1^+(2|2) \backslash{\Gamma_{\Lambda_{\overline{\eta}; \mu}}}}(1 + e^{-\beta})\right)
\end{align*}
\begin{align*}
 &= \frac{1}{e^{\epsilon_1}-e^{\epsilon_1}}. \frac{1}{e^{\delta_1} - e^{\delta_1}}
\left[ \sum_{w \in S_{2}\times S_2} \epsilon(w) w\left(e^{(2,0;1,-1) } \prod_{\beta\in\Delta_1^+(2|2) \backslash{\Gamma_{\Lambda_{\overline{\eta}; \mu}}}}(1 + e^{-\beta})\right)\right.\\
&\left. -\frac{1}{2}\sum_{w \in S_{2}\times S_2} \epsilon(w) w\left(e^{(0,0;1,1) } \prod_{\beta\in\Delta_1^+  (2|2) \backslash{\Gamma_{\Lambda_{\overline{\eta}; \mu}}}}(1 + e^{-\beta})\right) \right].
\end{align*}

Hence
 \begin{align*}
 \chi_{\underline{r}}({\rm ch} V\; (\Lambda_{\overline\eta; \mu}))&= \frac{1}{x_{r_1}-x_{r_2}}. \frac{1}{y_1 - y_2}
\left[ \chi_{\underline{r}}\left(  \sum_{w \in S_{2}\times S_2} \epsilon(w) w\left(e^{(2,0;1,-1) } \prod_{\beta\in\Delta_1^+(2|2) \backslash{\Gamma_{\Lambda_{\overline{\eta}; \mu}}}}(1 + e^{-\beta})\right)\right. \right)\\
&\left. -\frac{1}{2} \chi_{\underline{r}}\left( \sum_{w \in S_{2}\times S_2} \epsilon(w) w\left(e^{(0,0;1,1) } \prod_{\beta\in\Delta_1^+  (2|2) \backslash{\Gamma_{\Lambda_{\overline{\eta}; \mu}}}}(1 + e^{-\beta})\right) \right) \right].
\end{align*}

We have
\begin{align*}
&\chi_{\underline{r}}\left(  \sum_{w \in S_{2}\times S_2} \epsilon(w) w\left(e^{(2,0;1,-1) } \prod_{\beta\in\Delta_1^+(2|2) \backslash{\Gamma_{\Lambda_{\overline{\eta}; \mu}}}}(1 + e^{-\beta})\right) \right)\\
&=\sum_{w \in S_{2}\times S_2} \epsilon(w) w\left(x_{r_1}x_{r_2}^{-1}y_{1}y_{2}^{-1}(x_{r_1}+y_1 )  (x_{r_2}+ y_2 )  \right)\\
& =\sum_{w_1 \in S_{2}} x_{r_1}x_{r_2}^{-1}y_{1}y_{2}^{-1}(x_{r_1}+ y_1 )  (x_{r_2}+ y_2 ) - \sum_{w_1 \in S_{2}}x_{r_2}x_{r_1}^{-1}y_{1}y_{2}^{-1}(x_{r_2}+y_1 )  (x_{r_1}+ y_2 )
\end{align*}
and
\begin{align*}
&\chi_{\underline{r}}\left( \sum_{w \in S_{2}\times S_2} \epsilon(w) w\left(e^{(0,0;1,1) } \prod_{\beta\in\Delta_1^+  (2|2) \backslash{\Gamma_{\Lambda_{\overline{\eta}; \mu}}}}(1 + e^{-\beta})\right) \right)\\
&=\sum_{w \in S_{2}\times S_2} \epsilon(w) w\left(x_{r_1}^{-1}x_{r_2}^{-1}y_{1}y_{2}(x_{r_1}+y_1 )  (x_{r_2}+ y_2 )  \right)\\
& = \sum_{w_1 \in S_{2}}x_{r_1}^{-1}x_{r_2}^{-1}y_{1}y_{2}(x_{r_1}+ y_1 )  (x_{r_2}+ y_2 )- \sum_{w_1 \in S_{2}}x_{r_2}^{-1}x_{r_1}^{-1}y_{1}y_{2}(x_{r_2}+y_1 )  (x_{r_1}+ y_2 ),
\end{align*}
where $w_1$ only acts on the indices of $y$.

The second, we Computate $ch\;V(\Lambda_{\overline\kappa})$. Fallows by Lemma\eqref{L:kapa}
\begin{align*}
ch\; V(\Lambda_{\bar{\kappa}}) &= \frac{\prod_{\beta\in\Delta_1^+}(1|2) (1 + e^{-\beta})}{\prod_{i < j\leq 1}(e^{(\epsilon_{ i}- \epsilon_{ j})/2}- e^{-(\epsilon_{ i}- \epsilon_{ j})/2})} \sum_{w \in S_1}\epsilon(w) w(e^{\Lambda_{\bar{\kappa}} + \frac{1}{2}(0; 0, 0)})\\ 
& = \frac{(e^{\epsilon_3}+ e^{\delta_1} )  (e^{\epsilon_3}+ e^{\delta_2} )} {e^{2\epsilon_3}} e^{\Lambda_{\bar{\kappa}} }
\end{align*}
Thus
$$ \chi_{\underline{s}}( ch\; V(\Lambda_{\bar{\kappa}}))=   x_{r_3}^{-3}y_1^{-1}y_2^{-1}(x_{r_3}+ y_1 )  (x_{r_3}+ y_2 ).$$

The right hand side of \eqref{Exa}
\begin{align*}
 &\sum_{\underline{r}}\frac{ x_{r_1} x_{r_2} \chi_{\underline{r}}(ch\;V(\Lambda_{\overline\eta;\mu}) ) \chi_{\underline{s}}(ch\;V(\Lambda_{\overline\kappa}))}{ (x_{r_1}- x_{r_{3}}) (x_{r_2}- x_{r_{3}})}
= \sum_{\underline{r}} \frac{x_{r_1}x_{r_2}}{(x_{r_1}-x_{r_3})(x_{r_2}-x_{r_3})}  \chi_{\underline{s}}(ch\;V(\Lambda_{\overline\kappa}))\\
&\left[ \frac{1}{x_{r_1}-x_{r_2}}. \frac{1}{y_1 - y_2} \left(\sum_{w_1 \in S_{2}} x_{r_1}x_{r_2}^{-1}y_{1}y_{2}^{-1}(x_{r_1}+ y_1 )  (x_{r_2}+ y_2 ) - \sum_{w_1 \in S_{2}}x_{r_2}x_{r_1}^{-1}y_{1}y_{2}^{-1}(x_{r_2}+y_1 )  (x_{r_1}+ y_2 )\right)\right.\\
&\left. -\frac{1}{2}\frac{1}{x_{r_1}-x_{r_2}}. \frac{1}{y_1 - y_2} \left(\sum_{w_1 \in S_{2}}x_{r_1}^{-1}x_{r_2}^{-1}y_{1}y_{2}(x_{r_1}+ y_1 )  (x_{r_2}+ y_2 )- \sum_{w_1 \in S_{2}}x_{r_2}^{-1}x_{r_1}^{-1}y_{1}y_{2}(x_{r_2}+y_1 )  (x_{r_1}+ y_2 ) \right)\right]\\
&= Rhs1 -\frac{1}{2}Rhs2,
\end{align*}
where
\begin{align*}
Rhs1&:= \sum_{\underline{r} }\frac{x_{r_1}x_{r_2}}{(x_{r_1}-x_{r_3})(x_{r_2}-x_{r_3})} \chi_{\underline{s}}( ch\;V(\Lambda_{\overline\kappa}))\\
&\frac{1}{x_{r_1}-x_{r_2}}. \frac{1}{y_1 - y_2} \left(\sum_{w_1 \in S_{2}} x_{r_1}x_{r_2}^{-1}y_{1}y_{2}^{-1}(x_{r_1}+ y_1 )  (x_{r_2}+ y_2 ) - \sum_{w_1 \in S_{2}}x_{r_2}x_{r_1}^{-1}y_{1}y_{2}^{-1}(x_{r_2}+y_1 )  (x_{r_1}+ y_2 )\right)
\end{align*}
and
\begin{align*}
Rhs2&:= \sum_{\underline{r} \subset \{1, 2,3\}}\frac{x_{r_1}x_{r_2}}{(x_{r_1}-x_{r_3})(x_{r_2}-x_{r_3})}   \chi_{\underline{s}}( ch\;V(\Lambda_{\overline\kappa}))\\
&\frac{1}{x_{r_1}-x_{r_2}}. \frac{1}{y_1 - y_2}\left(\sum_{w_1 \in S_{2}}x_{r_1}^{-1}x_{r_2}^{-1}y_{1}y_{2}(x_{r_1}+ y_1 )  (x_{r_2}+ y_2 )- \sum_{w_1 \in S_{2}}x_{r_2}^{-1}x_{r_1}^{-1}y_{1}y_{2}(x_{r_2}+y_1 )  (x_{r_1}+ y_2 ) \right).
\end{align*}
We have
\begin{align*}
Rhs1&=\sum_{\underline{r} }\frac{1}{(x_{r_1}-x_{r_2})(x_{r_1}-x_{r_3})(x_{r_2}-x_{r_3})(y_1 - y_2)}   \chi_{\underline{s}}( ch\;V(\Lambda_{\overline\kappa}))\\
& \left( \sum_{w_1 \in S_{2}}x_{r_1}^2y_{1}y_{2}^{-1}(x_{r_1}+ y_1 )  (x_{r_2}+ y_2 )  - \sum_{w_1 \in S_{2}}x_{r_2}^2y_{1}y_{2}^{-1}(x_{r_2}+y_1 )  (x_{r_1}+ y_2 ) \right).
\end{align*}

Expand this right hand side, we have
\begin{align*}
Rhs1&=\sum_{\underline{r} }\frac{1}{(x_{r_1}-x_{r_2})(x_{r_1}-x_{r_3})(x_{r_2}-x_{r_3})(y_1 - y_2)} \\
&\left( \left(\sum_{w_1 \in S_{2}}x_{r_1}^2x_{r_3}^{-3}y_{2}^{-2}(x_{r_1}+ y_1 )  (x_{r_2}+ y_2 ) - \sum_{w_1 \in S_{2}}x_{r_2}^2x_{r_3}^{-3}y_{2}^{-2}(x_{r_2}+y_1 )  (x_{r_1}+ y_2 ) \right)(x_{r_3}+ y_1 )  (x_{r_3}+ y_2 ) \right)
\end{align*}
Since $\underline{r} \subset \{1, 2,3\}$ then  $\underline{r}$ equal to $\{1, 2\}or \{1, 3\}$ or $\{ 2,3 \}$. By replacing $\underline{r}$ by $\{r_1=1, r_2 =2\} $ or $\{r_1 =1, r_2 =3 \}$  or $ \{ r_1=3, r_2=2 \}$ we obtain\\
\begin{align*}
Rhs1&=\frac{1}{(x_{1}-x_{2})(x_{1}-x_{3})(x_{2}-x_{3})(y_1 - y_2)} \\
&\left( \left(\sum_{w_1 \in S_{2}}x_{1}^2x_{3}^{-3}y_{2}^{-2}(x_{1}+ y_1 )  (x_{2}+ y_2 ) - \sum_{w_1 \in S_{2}}x_{2}^2x_{3}^{-3}y_{2}^{-2}(x_{2}+y_1 )  (x_{1}+ y_2 ) \right)(x_{3}+ y_1 )  (x_{3}+ y_2 ) \right)
\end{align*}
\begin{align*}
&+ \frac{1}{(x_{1}-x_{3})(x_{1}-x_{2})(x_{3}-x_{2})(y_1 - y_2)} \\
&\left( \left(\sum_{w_1 \in S_{2}}x_{1}^2x_{2}^{-3}y_{2}^{-2}(x_{1}+ y_1 )  (x_{3}+ y_2 ) - \sum_{w_1 \in S_{2}}x_{3}^2x_{2}^{-3}y_{2}^{-2}(x_{3}+y_1 )  (x_{1}+ y_2 ) \right)(x_{2}+ y_1 )  (x_{2}+ y_2 ) \right)
\end{align*}
\begin{align*}
&+\frac{1}{(x_{3}-x_{2})(x_{3}-x_{1})(x_{2}-x_{1})(y_1 - y_2)} \\
&\left( \left(\sum_{w_1 \in S_{2}}x_{3}^2x_{1}^{-3}y_{2}^{-2}(x_{3}+ y_1 )  (x_{2}+ y_2 )  - \sum_{w_1 \in S_{2}}x_{2}^2x_{1}^{-3}y_{2}^{-2}(x_{2}+y_1 )  (x_{3}+ y_2 )  \right)(x_{1}+ y_1 )  (x_{1}+ y_2 ) \right).
\end{align*}

Now, we rewrite
\begin{align*}
Rhs1&= \frac{1}{\prod_{1\leq i<j \leq 3}(x_i-x_j)}\frac{1}{y_1-y_2} \\
&\left( \left(\sum_{w_1 \in S_{2}}x_{1}^2x_{3}^{-3}y_{2}^{-2}(x_{1}+ y_1 )  (x_{2}+ y_2 ) - \sum_{w_1 \in S_{2}}x_{2}^2x_{3}^{-3}y_{2}^{-2}(x_{2}+y_1 )  (x_{1}+ y_2 ) \right)(x_{3}+ y_1 )  (x_{3}+ y_2 ) \right)
\end{align*}
\begin{align*}
&+ \frac{1}{\prod_{1\leq i<j \leq 3}(x_i-x_j)}\frac{1}{y_1-y_2} \\
&\left( \left( -\sum_{w_1 \in S_{2}}x_{1}^2x_{2}^{-3}y_{2}^{-2}(x_{1}+ y_1 )  (x_{3}+ y_2 ) +\sum_{w_1 \in S_{2}} x_{3}^2x_{2}^{-3}y_{2}^{-2}(x_{3}+y_1 )  (x_{1}+ y_2 ) \right)(x_{2}+ y_1 )  (x_{2}+ y_2 ) \right)
\end{align*}
\begin{align*}
&+\frac{1}{\prod_{1\leq i<j \leq 3}(x_i-x_j)}\frac{1}{y_1-y_2} \\
&\left( \left( - \sum_{w_1 \in S_{2}}x_{3}^2x_{1}^{-3}y_{2}^{-2}(x_{3}+ y_1 )  (x_{2}+ y_2 )  +  \sum_{w_1 \in S_{2}}x_{2}^2x_{1}^{-3}y_{2}^{-2}(x_{2}+y_1 )  (x_{3}+ y_2 )  \right)(x_{1}+ y_1 )  (x_{1}+ y_2 ) \right).
\end{align*}

In the orther words
\begin{align*}
Rhs1&= \frac{1}{\prod_{1\leq i<j \leq 3}(x_i-x_j)}\frac{1}{y_1-y_2} \\
&\left[  w^1\left(\sum_{w_1 \in S_{2}}x_{1}^2x_{3}^{-3}y_{2}^{-2}(x_{1}+ y_1 )  (x_{2}+ y_2 )(x_{3}+ y_1 )  (x_{3}+ y_2 )\right) \right.\\
 &- w^2 \left(\sum_{w_1 \in S_{2}}x_{1}^2x_{3}^{-3}y_{2}^{-2}(x_{1}+ y_1 )  (x_{2}+ y_2 )(x_{3}+ y_1 )  (x_{3}+ y_2 )\right) 
\end{align*}
\begin{align*}
&+ w^3\left(\sum_{w_1 \in S_{2}}x_{1}^2x_{3}^{-3}y_{2}^{-2}(x_{1}+ y_1 )  (x_{2}+ y_2 )(x_{3}+ y_1 )  (x_{3}+ y_2 )\right) \\
 &-  w^4 \left(\sum_{w_1 \in S_{2}}x_{1}^2x_{3}^{-3}y_{2}^{-2}(x_{1}+ y_1 )  (x_{2}+ y_2 )(x_{3}+ y_1 )  (x_{3}+ y_2 )\right)
\end{align*}
\begin{align*}
&  w^5\left(\sum_{w_1 \in S_{2}}x_{1}^2x_{3}^{-3}y_{2}^{-2}(x_{1}+ y_1 )  (x_{2}+ y_2 )(x_{3}+ y_1 )  (x_{3}+ y_2 )\right) \\
 &\left.-  w^6\left (\sum_{w_1 \in S_{2}}x_{1}^2x_{3}^{-3}y_{2}^{-2}(x_{1}+ y_1 )  (x_{2}+ y_2 )(x_{3}+ y_1 )  (x_{3}+ y_2 )\right) \right],
\end{align*}
where $w^1 = (1)(1)(1)=(1), w^2= (12)(1)(1)=(12); w^3 =(1) (1)(23), w^4=(13)(1)(23) = (132); w^5=(1)(1) (13)=(13),w^6 =(23)(1)(13)=(123) $ are the permutation of $\{1,2,3\}$ and $w^i$ only acts on the indices of $x$.
From this we have
\begin{align*}
Rhs1&=\frac{1}{\prod_{1\leq i<j \leq 3}(x_i-x_j)}\frac{1}{y_1-y_2}\left[\sum_{w_0 \in S_3}\epsilon (w_0)w_0 \left(\sum_{w_1 \in S_{2}} x_1^2 x_3^{-3} y_2^{-2}(x_1+y_1)(x_2+y_2)(x_3+y_1)(x_3+y_2) \right)\right]\\
&=\frac{1}{\prod_{1\leq i<j \leq 3}(x_i-x_j)}\frac{1}{y_1-y_2}\left[\sum_{w \in S_3\times S_2}\epsilon (w)w \left( x_1^2 x_3^{-3} y_2^{-2}(x_1+y_1)(x_2+y_2)(x_3+y_1)(x_3+y_2) \right)\right].
\end{align*}

And
\begin{align*}
Rhs2&:= \sum_{\underline{r} \subset \{1, 2,3\}}\frac{x_{r_1}x_{r_2}}{(x_{r_1}-x_{r_3})(x_{r_2}-x_{r_3})}   \chi_{\underline{s}}( ch\;V(\Lambda_{\overline\kappa}))\\
&\frac{1}{x_{r_1}-x_{r_2}}. \frac{1}{y_1 - y_2}\left(\sum_{w_1 \in S_{2}}x_{r_1}^{-1}x_{r_2}^{-1}y_{1}y_{2}(x_{r_1}+ y_1 )  (x_{r_2}+ y_2 )- \sum_{w_1 \in S_{2}}x_{r_2}^{-1}x_{r_1}^{-1}y_{1}y_{2}(x_{r_2}+y_1 )  (x_{r_1}+ y_2 ) \right).
\end{align*}

So
\begin{align*}
Rhs2&=\sum_{\underline{r} \subset \{1, 2,3\}}\frac{1}{(x_{r_1}-x_{r_2})(x_{r_1}-x_{r_3})(x_{r_2}-x_{r_3})(y_1 - y_2)}   \chi_{\underline{s}}( ch\;V(\Lambda_{\overline\kappa}))\\
&  \left(\sum_{w_1 \in S_{2}}y_{1}y_{2}(x_{r_1}+ y_1 )  (x_{r_2}+ y_2 )  - \sum_{w_1 \in S_{2}}y_{1}y_{2}(x_{r_2}+y_1 )  (x_{r_1}+ y_2 ) \right).
\end{align*}
The calculation is similar to the case of Rhs1, we have
\begin{align*}
Rhs2&=\frac{1}{\prod_{1\leq i<j \leq 3}(x_i-x_j)}\frac{1}{y_1-y_2} \left(\sum_{w \in S_3 \times S_2}\epsilon (w)w \left(  x_3^{-3}(x_1+y_1)(x_2+y_2)(x_3+y_1)(x_3+y_2) \right)\right).
\end{align*}

The last, we compute ${\rm ch} V\; (\Lambda)$. We have
\begin{align*}
ch V(\Lambda) = \sum_{\pi \in S_2}\frac{1}{2!}\binom{2}{\pi}&(-1)^{|\Lambda - (\pi . \Lambda)_{\uparrow}| + l(\pi)} \frac{1}{L_0}\\
&\sum_{w \in S_3\times S_2}\epsilon (w)w(e^{(\pi . \Lambda)_{\uparrow} + \rho_0(3|2)}
\prod_{\beta \in \Delta_1^+ \backslash{\Gamma_{\Lambda}}}(1 + e^{-\beta}))
\end{align*}

 \begin{align*}
 &= \frac{1}{\prod_{1\leq i<j \leq 3}(e^{\epsilon_i}-e^{\epsilon_j})}\frac{1}{e^{\delta_1}-e^{\delta_2}} \left[\sum_{w \in S_3 \times S_2}\epsilon (w)w \left(e^{(3,1,-1;0,-2)}\prod_{\beta \in \Delta_1^+ \backslash{\Gamma_{\Lambda}}}(1 + e^{-\beta})\right)\right. \\
 &\left. -\frac{1}{2}\sum_{w \in S_3 \times S_2}\epsilon (w)w \left(e^{(1,1,-1;0,0)}\prod_{\beta \in \Delta_1^+ \backslash{\Gamma_{\Lambda}}}(1 + e^{-\beta})\right)\right].
\end{align*}

Consider
\begin{align*}
 &\frac{1}{\prod_{1\leq i<j \leq 3}(e^{\epsilon_i}-e^{\epsilon_j})}\frac{1}{e^{\delta_1}-e^{\delta_2}}\sum_{w \in S_3 \times S_2}\epsilon (w)w \left(e^{(3,1,-1;0,-2)}\prod_{\beta \in \Delta_1^+ \backslash{\Gamma_{\Lambda}}}(1 + e^{-\beta})\right) \\
&=\frac{1}{\prod_{1\leq i<j \leq 3}(x_i-x_j)}\frac{1}{y_1-y_2} \sum_{w \in S_3 \times S_2}\epsilon (w)w \left( x_1^3x_2x_3^{-1}y_2^{-2} \frac{(x_1+y_1)(x_2+y_2)(x_3+y_1)(x_3+y_2)}{x_1x_2x_3^2} \right) \\
 &=\frac{1}{\prod_{1\leq i<j \leq 3}(x_i-x_j)}\frac{1}{y_1-y_2}\sum_{w \in S_3 \times S_2}\epsilon (w)w \left( x_1^2 x_3^{-3} y_2^{-2}(x_1+y_1)(x_2+y_2)(x_3+y_1)(x_3+y_2) \right)\\
&=Rhs1.
\end{align*}

And
\begin{align*}
 &\frac{1}{\prod_{1\leq i<j \leq 3}(e^{\epsilon_i}-e^{\epsilon_j})} \frac{1}{e^{\delta_1}-e^{\delta_2}} \sum_{w \in S_3 \times S_2}\epsilon (w)w \left(e^{(1,1,-1;0,0)}\prod_{\beta \in \Delta_1^+ \backslash{\Gamma_{\Lambda}}}(1 + e^{-\beta})\right).\\
&=\frac{1}{\prod_{1\leq i<j \leq 3}(x_i-x_j)}\frac{1}{y_1-y_2} \sum_{w \in S_3 \times S_2}\epsilon (w)w \left( x_1x_2x_3^{-1} \frac{(x_1+y_1)(x_2+y_2)(x_3+y_1)(x_3+y_2)}{x_1x_2x_3^2} \right) \\
&=\frac{1}{\prod_{1\leq i<j \leq 3}(x_i-x_j)}\frac{1}{y_1-y_2}\sum_{w \in S_3 \times S_2}\epsilon (w)w \left( x_3^{-3} (x_1+y_1)(x_2+y_2)(x_3+y_1)(x_3+y_2) \right) \\
&= Rhs2.
\end{align*}
Thus
$$  ch V(\Lambda)=Rhs1 - \frac{1}{2}Rhs2.$$
Finally, we have
$$  ch V(\Lambda)= \sum_{\underline{r}}\frac{ x_{r_1} x_{r_2} \chi_{\underline{r}}(ch\;V(\Lambda_{\overline\eta;\mu}) ) \chi_{\underline{s}}(ch\;V(\Lambda_{\overline\kappa}))}{ (x_{r_1}- x_{r_{3}}) (x_{r_2}- x_{r_{3}})}.$$

\end{proof}

\section{Supersymmetric S-fuctions}

\subsection{Symmetric functions associated to composite partitions}
Let $x = (x_1, \ldots, x_m)$ be a set of variables. Let $\bar{\nu}; \mu$ be a composite partition with lengths $l(\mu)= p, l(\nu) = q$ such that $p + q \leq m$. Thus $\bar{\nu}; \mu$ is an $m$-standard composite partition (cf. \eqref{eq:42}). We can associate to it an $m$-tuple $(\mu_1, \ldots, \mu_p,0, \ldots, 0, -\nu_q, \ldots, -\nu_1)$ with the composie partition $\bar{\nu}; \mu$. The symmetric Schur function indexed by this composite partition is defined by:
\begin{equation}\label{eq:51}
s_{\bar{\nu}; \mu}(x) = (\prod_{i=1}^mx_i^{-\nu_1})s_{\lambda}(x) ,
\end{equation} 
where $\lambda$ is the partition of length $m$ defined by $$ (\mu_1 + \nu_1, \mu_2 + \nu_1, \ldots, \mu_p+\nu_1, \nu_1,\ldots,\nu_1,-\nu_q +\nu_1,\ldots,-\nu_2 + \nu_1, 0 ).$$

A formula for symmetric Schur function indexed by the composite partition was postulated by Balantekin and Bars \cite{bar2} in terms of characters, and proved in \cite{CK}, namely
\begin{equation}\label{eq:52}
s_{\bar{\nu};\mu} (x)= det\left(
\begin{tabular}{c|c}
$\dot{h}_{\nu_l+k-l}(x)$&$h_{\mu_j-k-j+1}(x)$\\
\hline
$\dot{h}_{\nu_l-i-l+1}(x)$&$h_{\mu_j +i -j}(x)$
\end{tabular}\right),
\end{equation}
where the indices $i, j, k$ resp. $l$ run from top to bottom, from left to right, from bottom to top resp. from right to left and  $\dot{h_r}(x) = h_r(\bar{x})=  h_r(x_1^{-1}, \ldots, x_m^{-1})$.\\

Let $x'$  and $ x''$  be two subsets of $\{x_1, x_2, \ldots, x_m\}$. Denote by $E(x', x'')=\prod_{x_i \in x'}\prod_{x_j\in x''}(x_i-x_j)$ and $|x'| $ (resp. $|x''| $) is size of $x'$ (resp. $ x''$ ).
\begin{lem}\label{L:61}
For $m = p+q$, let $\mu= (\mu_1, \mu_2, \ldots, \mu_p), \nu= (\nu_1, \nu_2, \ldots, \nu_q)$ be two partitions and $\lambda= (\mu_1, \ldots, \mu_p, \nu_1, \ldots, \nu_q)$. Then 
\begin{equation}
\sum_{x' ,x''}\frac{s_{ \mu +q^p}(x') s_{\nu }( x'')}{E(x', x'')} = s_{\lambda}(x),
\end{equation}
where the sum is over all possible decompositions $\{x_1, x_2, \ldots, x_m \} = x' \cup x'' $ with $|x'| = p $, $|x''| = q$.
\end{lem}
\begin{proof}

This follows from \cite[Lemma 3.3]{M1}.
\end{proof}

\subsection{Super-symmetric functions associated to composite partitions}
Let $x = (x_1, x_2, \ldots, x_m) $ and $ y = y^{(n) }= (y_1, y_2, \ldots, y_n)$ be two sets of independent variables. The \textit{\textit{complete supersymmetric} functions }can be expressed in terms of the elementery symmetric and the complete symmetric functions:
\begin{equation}\label{eq:}
h_r(x/y) = \sum_{k = 0}^r h_k(x)e_{r-k}(y).
\end{equation} 

Given the complete supersymmetric functions $h_r(x/y)$ and any partition $\lambda = (\lambda_1, \lambda_2, \ldots)$, one defines the \textit{supersymmetric Schur fuctions} $s_\lambda(x/y)$:
\begin{equation}\label{eq:53}
s_\lambda(x/y) = det (h_{\lambda_i - i + j}(x/y))_{1\leq i,j \leq l(\lambda)}.
\end{equation} 
Obviously,
\begin{equation}\label{eq:}
s_{(r)}(x/y) = h_r(x/y).
\end{equation} 

Given a composite partition $\bar{\nu}; \mu$, one can define the corresponding \textit{supersymmetric Schur function}, also called \textit{supersymmetric $S$-function}\cite{bar1,bar2}:\\
\begin{equation}\label{eq:54}
 s_{\bar{\nu};\mu} (x/y)= det\left(
\begin{tabular}{c|c}
$\dot{h}_{\nu_l+k-l}(x/y)$&$h_{\mu_j-k-j+1}(x/y)$\\
\hline
$\dot{h}_{\nu_l-i-l+1}(x/y)$&$h_{\mu_j +i -j}(x/y)$
\end{tabular}\right),
\end{equation}
where the indices $i, j, k$ resp. $l$ run from top to bottom, from left to right, from bottom to top resp. from right to left and $\dot{h_r}(x/y) = h_r(\bar{x}/\bar{y})$ with $\bar{x}=(x_1^{-1}, \ldots, x_m^{-1}), \bar{y}=(y_1^{-1}, \ldots, y_n^{-1})$. \\
 For $ \nu = 0$, this supersymmetric $S$-function is so-called supersymmetric Schur function as  defined in \eqref{eq:53}.

\begin{lem}\label{L:62}
 Let $\bar{\nu}; \mu$ be a composite partition. Then
$$ s_{\bar{\nu}; \mu}(x/y) = \sum_{\alpha, \beta}s_{\bar{\beta}; \alpha}(x/y^{(n-1)}) y_n^{a-b},$$
where $a= \left|\mu -\alpha \right|, b= \left|\nu -\beta \right|$ and the sum is taken over all partitions $\alpha$ and $\beta$ such that $ (\mu - \alpha)_i, (\nu - \beta)_i \in \{0,1\}$.
\end{lem}
\begin{proof}
This follows from \cite[Lemma A.3]{M4}.
\end{proof}
\subsection{A recurrent formula for supersymmetric $S$-function}
\begin{lem}
Let $\alpha, \mu$ be partitions such that $\alpha \subset \mu $. Then
$$ m  - l(\mu) \leq m  - l(\alpha),$$
for all positive integer $m$. 
\end{lem}
\begin{proof}[Proof]
Because $\alpha \subset \mu$ so $ l(\alpha) \leq l(\mu)$. Thus
$$  m - l(\mu) \leq  m - l(\alpha).$$
\end{proof}
  
 \begin{thm}\label{T:S-function}
 Let $\bar{\nu};\mu$ be a composite partition. Assume $q$ is a nonnegative integer such that $ 0 < q < m + 1 - l(\mu) $ and $p = m - q$. Then
\begin{equation}\label{eq:S-function}
\sum_{x' ,x''}\frac{(\prod x')^q s_{\overline\eta; \mu}(x'/y) s_{\overline\kappa }( x''/ y)}{E(x', x'')} = s_{\overline \nu; \mu}(x/y),
\end{equation}
where $ \kappa = (\nu_1, \nu_2, \ldots, \nu_q)$, $\eta = (\nu_{q + 1}, \nu_{q + 2}, \ldots)$ and the sum is over all possible decompositions $\{x_1, x_2, \ldots, x_m\} = x' \cup x'' $ with $|x'| = p $ , $|x''| = q$.
\end{thm}

\begin{proof}[Proof] We use induction on $n$.  Let's consider the case $ n = 0 $. We consider the following subcases:\\
Subcase 1: $l(\mu) + l(\nu) - m -1 < 0$, in the orther words $\overline\nu; \mu$ is $m$-standard composite partition. Since $ q < m + 1 - l(\mu) $ hence $\overline\eta; \mu$ is $p$-standard composite partition. The left hand side of \eqref{eq:S-function} is:
\begin{align*}
 Lhs: & = \sum_{x' ,x''}\frac{(\prod x')^q s_{\overline\eta; \mu }(x') s_{\overline\kappa }( x'')}{E(x', x'')}\\
&= \sum_{x' ,x''}\frac{(\prod x')^q (\prod x')^{-\eta_1} s_{(\mu_1 + \eta_1, \mu_2 + \eta_1, \ldots, -\eta_2 + \eta_1, 0) }(x') s_{\overline\kappa}( x'')}{E(x', x'')}\\
&= (\prod x)^{- \nu_1}\sum_{x' ,x''}\frac{ (\prod x')^q (\prod x')^{-\eta_1} (\prod x)^{ \nu_1} s_{(\mu_1 + \eta_1, \mu_2 + \eta_1, \ldots, -\eta_2 + \eta_1, 0) }(x')  s_{\overline\kappa}( x'')}{E(x', x'')}\\
&= (\prod x)^{- \nu_1}\sum_{x' ,x''}\frac{ s_{(\mu_1 + \nu_1, \mu_2 + \nu_1, \ldots, -\eta_2 + \nu_1, -\eta_1 + \nu_1) + (q^p)}(x')s_{\overline\kappa + (\nu^q_1)}( x'')}{E(x', x'')}\\
&= (\prod x)^{- \nu_1}s_{(\mu_1  + \nu_1, \mu_2  + \nu_1, \ldots, -\nu_2 + \nu_1, 0)}(x),
 \end{align*}
where the last equality follows  Lemma \ref{L:61}. We deduce\\
 
                 $$ \sum_{x' ,x''}\frac{(\prod x')^q s_{\overline\eta; \mu }(x')s_{\overline\kappa }( x'')}{E(x', x'')} = s_{\overline \nu; \mu}(x). $$\\

Subcase 2: $ 0 \leq l(\mu) + l(\nu) - m -1 < l(\nu)$. We put  
$$q_0 = l(\mu) + l(\nu) - m -1,$$
then
$$ 0\leq q_0 =  l(\mu) + l(\nu) - m -1 = l(\nu) - ( m + 1 - l(\mu)) \leq l(\nu) - q = l(\eta).$$
On the orther hand

\begin{align*}
q_0&= l(\mu) + l(\nu) - m - 1\\ 
&= l(\mu) + q + l(\eta) - m -1 \\
&= l(\mu) + l(\eta) - (m - q) - 1\\
&= l(\mu) + l(\eta) - p - 1.
\end{align*} 
By applying modification rule (see \cite{CK}), we have

$$
\begin{cases}
     s_{\overline\nu;\mu}(x) = (-1)^{c + \overline c + 1}s_{\overline{\nu - q_0}; \mu - q_0}(x)\\
s_{\overline\eta;\mu }(x) =  (-1)^{c + \overline c + 1}s_{\overline{\eta - q_0}; \mu  - q_0}(x).
\end{cases} 
$$
Here $ \nu - q_0, \mu - q_0 $ ( resp. $ \eta - q_0$) are derived from Young diagrams $F^{\nu}, F^{\mu}$  (resp. $F^{ \eta} $) by removing $ q_0 $ cells which are continuous boundaries starting at the foot of the first column of the Young diagram $F^{\nu}, F^{\mu}$ (resp. $F^{ \eta}$),  and ending at the $ \overline c $, $ c $ (resp. $ \overline c $ ) column of Young diagram $F^{\nu}, F^{\mu}$ (resp. $F^{ \eta}$).\\
By applying repeatedly the  modification rule  after a finite number of steps, say affter $k$ steps, we obtain:\\

$$
\begin{cases}
     s_{\overline\nu;\mu} (x)= \pm s_{\overline{ \beta};\alpha}(x)\\
s_{\overline\eta;\mu }(x) =  \pm s_{\overline{ \gamma}; \alpha )}(x),
\end{cases} 
$$
where $l(\alpha) + l(\beta) < m + 1, l(\alpha) + l(\gamma) < p + 1 $.\\

Here the number of times modification rule used to compute $ s_{\overline \nu; \mu}(x) $ is equal to the number of times that the used modification rule to compute $ s_{\overline \eta; \mu} (x) $. At each step, their two sign are the same. In addition, the numbers of elements removed in each column of $ \nu $ and $ \eta $ are the same. One have

$$
\begin{cases}
    \alpha = (\mu_1, \mu_2, \ldots,  \mu_{l(\mu) - (q_0 - 1)} - s_{(q_0 - 1)},  \mu_{l(\mu) - (q_0 - 2)} - s_{(q_0 - 2)}, \dots,\mu_{l(\mu) } - s_{0} ) \\
\beta = (\nu_1, \nu_2, \ldots, \nu_{l(\nu) - (q_0 - 1)} - t_{(q_0 - 1)},  \nu_{l(\nu) - (q_0 - 2)} - t_{(q_0 - 2)}, \dots,\nu_{l(\nu) } - t_{0} ) \\
\gamma = (\nu_{q +1}, \nu_{q +2} \ldots, \nu_{l(\nu) - (q_0 - 1)} - t_{(q_0 - 1)},  \nu_{l(\nu) - (q_0 - 2)} - t_{(q_0 - 2)}, \dots,\nu_{l(\nu) } - t_{0} ) 
\end{cases} 
$$
where $s_i$ and $t_i$, $i = 0, 1, \ldots, q_0 - 1$, are positive integers.\\

We are now reduced to subcase 1. Therefore, we have
\begin{equation*}
\sum_{x' ,x''} \frac{(\prod x')^q s_{\overline\gamma; \alpha}(x'/y)s_{\overline\kappa }( x''/ y)}{E(x', x'')} = s_{\overline \beta; \alpha}(x/y),
\end{equation*}
where $\kappa = (\beta_1, \beta_2, \ldots, \beta_q), \gamma = (\beta_{q + 1}, \beta_{q + 2}, \ldots)$.\\
So, we obtain

\begin{equation*}
\sum_{x' ,x''}\frac{(\prod x')^q s_{\overline\eta; \mu}(x'/y)s_{\overline\kappa }( x''/ y)}{E(x', x'')} = s_{\overline \nu; \mu}(x/y).
\end{equation*}
\\
Let consider the induction step. Supose the equality holds for all values $k < n$.  Denote $y^{(n)} = (y_1, y_2, \ldots, y_n)$. We use Lemma \ref{L:62} to isolate $y_n$, this gives 
\begin{align*}
Lhs: &=  \sum_{x' ,x''}\frac{(\prod x')^q s_{\overline\eta; \mu}(x'/y)s_{\overline\kappa }( x''/ y)}{E(x', x'')}    \\ 
& =  \sum_{x' ,x''}\frac{(\prod x')^q }{E(x', x'')} \left(\sum_{\alpha, \beta} s_{\overline\beta;\alpha}(x'/y^{(n-1)})y_n^{a - b}\right) \left(\sum_{\gamma} s_{\overline\gamma }( x''/ y^{(n-1)})y_n^{-c}\right),
\end{align*}
where $(\mu - \alpha)_i, (\eta - \beta)_i, (\kappa - \gamma)_i \in \{0, 1\}, |\mu - \alpha| = a, |\eta - \beta| = b, |\kappa - \gamma| = c$.
\begin{align*}
 &= \sum_{\alpha, \beta, \gamma} \left(\sum_{x' + x''}\frac{(\prod x')^q  }{E(x', x'')} s_{\overline\beta;\alpha}(x'/y^{(n-1)}) s_{\overline\gamma}( x''/ y^{(n-1)}) \right) y_n^{a - b - c}\\
 &= \sum_{\alpha, \tau} s_{\overline \tau; \alpha}(x/y^{(n- 1)}) y_n^{a - b - c} 
\end{align*}
where $|\mu - \alpha| = a, |\nu - \tau| = b + c, (\mu - \alpha)_i \in \{0,1\}, (\nu - \tau)_i \in \{0,1\}$. Further more, the last equality comes from the inductive hypothesis with composite partition $\overline \tau; \alpha$ satifies $q <  m + 1 - l(\mu) \leq m + 1 - l(\alpha)$. \\
Now, by reusing Lemma \ref{L:62} we have
 $$
Lhs= s_{\overline \nu; \mu}(x/y). $$
From this the statment follows for $n$.\\
\end{proof}

\section{Proof the main theorem}
Now we restate the main theorem (Theorem \ref{Mt}) and prove it.
\begin{thm} 
Let  $ \Lambda = (\alpha_1, \alpha_2, \ldots, \alpha_m; -k, -k, \ldots,  -k)$ be an integral dominant weight of a general linear Lie superalgebra $\mathfrak{gl}(m|n)$ such that $\alpha_{m-k} \geq 0 \geq \alpha_{m-k+1}$ with $k$ is a nonnegative integer and $0 \leq k \leq m$. Then
$$  {\rm ch} V(\Lambda) = det\left(
\begin{tabular}{c|c}
$\dot{h}_{n-\alpha_{m-t+1}+s-t}(x/y)$&$h_{\alpha_j-s-j+1}(x/y)$\\
\hline
$\dot{h}_{n-\alpha_{m-t+1}-i-t+1}(x/y)$&$h_{\alpha_j +i -j}(x/y)$
\end{tabular}\right),$$
where the indices $i, j, s$ resp. $t$ run from top to bottom, from left to right, from bottom to top resp. from right to left with $i,j = 1,2,\ldots, m-k$ and $s,t=1,2,\ldots, k$.
\end{thm}
 \begin{proof} For $ \Lambda = (\alpha_1, \alpha_2, \ldots, \alpha_m; -k, -k, \ldots,  -k)$ with $0 \leq k \leq m$ and $\alpha_{m-k}\geq 0 \geq \alpha_{m-k+1}$. Then $\Lambda$ is a special weight in $P_k$ and the $(m|n)$-standard composite partiton $\bar{\nu};\mu \in Q_k$ corresponds to $\Lambda$ is defined by 
$$
\begin{cases}
    \mu &= (\alpha_1, \alpha_2, \ldots, \alpha_{m - k}),   \\ 
    \nu &= (n-\alpha_m, n-\alpha_{m-1}, \ldots, n-\alpha_{m-k+1}), \\
\end{cases} 
$$
(cf. Proposition \ref{Pro:2} ).

By applying the Theorem \ref{T:ch} with $\eta=(0), \kappa = \nu$ we have
\begin{align}\label{eq:71}
 {\rm ch} V(\Lambda)  = \sum_{ \underline{r}, \underline{s}} \frac{(\prod_{i = 1}^{m-k}  x_{r_i})^k \chi_{\underline{r}}(ch\; V(\Lambda_{\overline0;\mu}))  \chi_{\underline{s}}(ch\;V(\Lambda_{\overline\nu}))}{\prod_{i= 1}^{m-k} \prod_{j = 1}^k (x_{r_i}- x_{r_{m-k+j}})}, 
\end{align}
where the sum is over all possible decomposition $\underline{r} \cup \underline{s}$ with $\underline{r} := \{r_1, r_2, \ldots, r_{m-k}\} \subset \{1, 2, \ldots, m\}$ and $ \underline{s}= \{r_{m-k + 1},r_{m-k + 2}, \ldots, r_{m}\} = \{1, 2, \ldots, m\} \backslash \underline{r}$.

By applying Theorem \ref{T:S-function} with $\bar{\nu};\mu$ as above, $\eta=0, \kappa= \nu$ and $q=k, p= m-k$ we have

\begin{align}\label{eq:72}
 \sum_{x' ,x''}\frac{(\prod x')^k s_{ \mu}(x'/y) s_{\overline\nu }( x''/ y)}{E(x', x'')} = s_{\overline \nu; \mu}(x/y),
\end{align}
 where the sum is over all possible decompositions $\{x_1, x_2, \ldots, x_m\} = x' \cup x'' $ with $|x'| = m-k $ , $|x''| = k$.

Now set $x'=\{x_{r_1},x_{r_2},\ldots,x_{r_{m-k}}\}$, $x''=\{x_{r_{m-k+1}},x_{r_{m-k+2}},\ldots,x_{r_{m}}\}$ then $\{x_1, x_2, \ldots, x_m\} = x' \cup x'' $ and as we know
$$ \chi_{\underline{r}}(ch\; V(\Lambda_{\overline0;\mu}))= s_{ \mu}(x'/y),
 $$
$$ \chi_{\underline{s}}(ch\;V(\Lambda_{\overline\nu}))= s_{\overline\nu }( x''/ y).
 $$
So
\begin{align*}
 \sum_{\underline{r}, \underline{s}} \frac{(\prod_{i = 1}^{m-k}  x_{r_i})^k  \chi_{\underline{r}}(ch\; V(\Lambda_{\mu}))  \chi_{\underline{s}}(ch\;V(\Lambda_{\overline\nu}))}{\prod_{i= 1}^{m-k} \prod_{j = 1}^k (x_{r_i}- x_{r_{m-k+j}})}
=\sum_{x' ,x''}\frac{(\prod x')^k s_{ \mu}(x'/y) s_{\overline\nu }( x''/ y)}{E(x', x'')}.
\end{align*}
From equation \eqref{eq:71} and equation \eqref{eq:72} deduce
$$  {\rm ch} V(\Lambda)  = s_{\overline \nu; \mu}(x/y).
 $$
Namely, in the formula \eqref{eq:SZ}, the formula used to calculate ${\rm ch} V(\Lambda)$ if we substitute $ e^{\delta_i}=x_i $ for $i=1,2,\ldots,m$ and $ e^{\delta_j}= y_j $ for $j= 1,2,\ldots,n$ then
$$  {\rm ch} V(\Lambda) = det\left(
\begin{tabular}{c|c}
$\dot{h}_{n-\alpha_{m-t+1}+s-t}(x/y)$&$h_{\alpha_j-s-j+1}(x/y)$\\
\hline
$\dot{h}_{n-\alpha_{m-t+1}-i-t+1}(x/y)$&$h_{\alpha_j +i -j}(x/y)$
\end{tabular}\right),$$
where the indices $i, j, s$ resp. $t$ run from top to bottom, from left to right, from bottom to top resp. from right to left with $i,j = 1,2,\ldots, m-k$ and $s,t=1,2,\ldots, k$.
\end{proof}

\section{Concluding remarks} 
The starting point of this work is the construction of irreducible representations of $\mathfrak{gl}(3|1)$, given in \cite{dhh,dung,dh}, which suggests a determinantal type formula expressing an irreducible representations in terms of the symmetric powers of the fundamental representation and their duals. We later discover that a vast generalization of this formula has been provided by Moens and van der Jeugt.
In the paper \cite[2004]{M4}, E.M. Moens and J. van der Jeugt announce the following theorem.\\

\noindent{\bf Theorem.} \cite[Theorem 4.3]{M4}  {\em
Let $\bar{\nu}; \mu$ be a standard and critical composite partition (see \cite[Definition 3.1]{M4}) with no overlap (see \cite[Section 2]{M4}) and $\Lambda_{\bar{\nu}; \mu}$ be the corresponding super weight. The character $  {\rm ch} V(\Lambda_{\bar{\nu}; \mu}) $ is equal to $s_{\bar{\nu};\mu}(x/y) $, which is defined in the same manner as in \eqref{eq:54}.
 }\\

The proof of this theorem is based on the following lemma. 

\noindent{\bf Lemma.}  \cite[Lemma A.5]{M4} {\em 
Suppose $|x| = m, |y| = n$ and $h, p$, $q$ are positive integers with $ m = p + q$. Let $\kappa=(\kappa_1, \kappa_2, \ldots, \kappa_q), \eta =(\eta_1, \eta_2, \ldots)$ and $\mu$ be partitions, and $\nu = (\kappa_1, \kappa_2, \ldots, \kappa_q, \eta_1, \eta_2, \ldots)$. Then
\begin{equation*}
\sum_{x' + x''}\frac{(\prod x')^q (\prod x'')^h s_{\overline\eta; \mu}(x'/y) s_{\kappa + (h^q)}( \bar{x''}/ \bar{y})}{E(x', x'')} = s_{\overline \nu; \mu}(x/y),
\end{equation*}
where the sum is over all possible decompositions $x = x' + x''$ with $|x'| = p $ , $|x''| = q$.
}\\

However, Moens notices in his thesis   that this Lemma is false and proposes to prove the above theorem by using a weaker form of this Lemma, in which one assumes that $\bar{\nu}; \mu$ is a critical composite partition with no zeros in the overlap when presented in the $m\times n$-rectangle \cite[Lemma 5.14]{M6}. 

However no proofs are provided.
 Thus the mentioned above theorem  in its general form is still  a conjecture.

\section{Acknowledgment}
The research of authors is supported in part by  Vietnam National Foundation for Science and Technology Development (NAFOSTED), grant number 101.04-2020.19.


\begin{thebibliography}{10}

\bibitem{br1}
A.~Berele and A.~Regev,
\newblock {Hook Young Diagrams with Applications to Combinatorics and to
  Representation of Lie Algebras}.
\newblock {\em Advances in Math.}, 64:118--175, 1987.

\bibitem{Brundan}
J. Brundan,
\newblock {Kazhdan-Lusztig polynomials and character formulae for the Lie superalgebra $\mathfrak{gl}(m|m)$},
\newblock {\em J. Amer. Math. Soc.}, 16:185-231, 2002.

\bibitem{bar1}
A.B. Balantekin and I. Bars,
\newblock {Dimension and character formulas for lie supergrous},
\newblock {\em J. Math. phys. 2}, 1149-1162, 1981.

\bibitem{bar2}
A.B. Balantekin and I. Bars,
\newblock {Representation of supergrous},
\newblock {\em J. Math. phys. 2}, 1810-1818, 1981.

\bibitem{CK}
C.J. Cummins and R.C. King,
\newblock {Composite Young diagrams, supercharacters of $U(M/N)$ and modification rules},
\newblock {\em J. Phys. A}, 1987, no.11, 3121-3133.

\bibitem{dj}
P.H. Dondi and P.D. Jarvis,
\newblock {Diagram and superfild techniques in the classical superalgebras},
\newblock {\em J. Phys. A14 (1981)}, 547-563.

\bibitem{bdh} N.L.T. Binh, N.T.P. Dung and P.H. Hai,
\newblock{\em Jacobi-Trudi type formula for character of  irreducible representations of  $\frak{gl}(m|1)$,}
\newblock{  Acta Math. Vietnam.}, 44 (2019), no. 3, 603–615.

\bibitem{dhh} N.T.P. Dung, P.H. Hai and N.H. Hung,
\newblock{\em Construction of irreducible representations of the quantum super group $GL_q(3|1)$,}
\newblock{  Acta Math. Vietnam.}, 36 (2011), no. 2, 215–229.

\bibitem{dung}  N.T.P. Dung,
\newblock{\em Double Koszul complex and construction of irreducible representations of $\mathfrak{gl}(3|1)$,}
\newblock{  Proc. Amer. Math. Soc.} 138 (2010), no. 11, 3783–3796. 

\bibitem{dh} N.T.P. Dung and  P.H. Hai,
\newblock{\em Irreducible representations of quantum linear groups of type $A_{1|0}$,}
\newblock{  J. Algebra} 282 (2004), no. 2, 809–830. 

\bibitem{H}
J.W.B Hughes, R.C. King and J. van der Juegt,
\newblock {On the coposition factors of Kac modules for the Lie supralgebra $\frak{sl(m|n)}$},
\newblock {\em J. Math. Phys. 33 (1992)}, 470-491.

\bibitem{Kac1}
V.G. Kac,
\newblock {Classification of simple Lie superalgebras},
\newblock {\em Funct.Anal. Appl.}, 9:263-265, 1975.


\bibitem{Kac2}
V.G. Kac,
\newblock { Lie superalgebras},
\newblock {\em Adv. Math.}, 26:8-96, 1977.

\bibitem{Kac3}
V.G. Kac,
\newblock {Character of typical  representations of classical Lie superalgebras},
\newblock {\em Comm. Alg.}, 5:889-897, 1977.

\bibitem{Kac4}
V.G. Kac,
\newblock {Representations of classical Lie superalgebras, }
\newblock {\em in: Lecture Notes in Math.}, 676:597-626, 1978.

\bibitem{Macdonald}
I.G.Macdonald,
\newblock { Symmetric Function and the Hall Polynomials},
\newblock {\em Oxford University Press, New York,1979}.


\bibitem{M1}
E.M. Moens and J. van der Jeugt,
\newblock { A detrminantal fomula for super-symmetric schur polynomials},
\newblock {J. Algebraic Combin.17(2003)}, no. 3, 283 – 307.

\bibitem{M2}
E.M. Moens and J. van der Jeugt,
\newblock {On dimension formulas for $\mathfrak{gl}(m|n)$ representations},
\newblock {J. Lie Theory 14 (2004),} , no. 2, 523 – 535.

\bibitem{M3}
E.M. Moens and J. van der Jeugt,
\newblock {On characters and dimension fomulas for representation of the Lie superalgebra $\mathfrak{gl}(m|n)$},
\newblock {Lie Theory and Its Applications in Physics V, ed. H.-D. Doebner and V.K. Dobrev, World Sci. Publ., Singapore (2004),}64 – 73.

\bibitem{M4}
E.M. Moens and J. van der Jeugt,
\newblock {A character formula for atypical critical $\mathfrak{gl}(m|n)$ representations labeled by composite partitions},
\newblock {J. Phys.A: Math. Gen. 37 (2004),} no. 2, 523 – 535.

\bibitem{M5}
E.M. Moens and J. van der Jeugt,
\newblock { Composite super-symmetric S-functions and character of $\mathfrak{gl}(m|n)$ representations},
\newblock {Proceedings of theVI International Worshop on Lie Theory and its Applications in Physics, ed. H.-D. Doebner and V.K. Dobrev, Heron Press Ltd, Sofia (2006)}, 251 – 268.

\bibitem{M6}
E.M. Moens,
\newblock {Supersymmetric Schur functions and Lie superalgebra representations},
\newblock {Ph.D. thesis, University of  Gent, 2006}.

%

\bibitem{Zhang2}
Y. Su, R.B. Zhang,
\newblock {Character and dimension formula for general linear superalgebra}.
\newblock {\em Adv. Math.}, 211:1-33, 2007.

\bibitem{Jeugt1}
J. van der Jeugt, J.W.B Hughes, R.C. King and J. Thierry-Mieg,
\newblock {Character fomulas for irreducible modules of the Lie superalgebra $\mathfrak{gl}(m|n)$}.
\newblock {\em J. Math. Phys}, 31, no.1, 2278-2304, 1990.


\bibitem{Jeugt2}
J. van der Jeugt, J.W.B Hughes, R.C. King and J. Thierry-Mieg,f
\newblock {Character fomulas for irreducible modules of the Lie superalgebra $\mathfrak{sl}(m|n)$}.
\newblock {\em J. Math. Phys}, 31, no.1, 2278-2304, 1990.
\end{thebibliography}
\end{document}